\documentclass[11pt]{article}
\usepackage[latin1]{inputenc}
\usepackage[english]{babel}
\usepackage{amssymb,amsfonts, amsmath, color}
\usepackage[margin=2.7cm]{geometry}
\usepackage{xcolor}
\usepackage{graphicx, color, enumerate}
\usepackage[latin1]{inputenc}
\usepackage[active]{srcltx}
\usepackage{tikz}
\usepackage{pgf}
\usepackage{etex}
\usepackage{verbatim}
\usepackage{tikz-3dplot}
\usepackage{pgfkeys}
\usepackage{fp}

\newtheorem{theorem}{Theorem}[section]
\newtheorem{proposition}[theorem]{Proposition}

\newtheorem{lemma}[theorem]{Lemma}
\newtheorem{remark}[theorem]{Remark}

\newcommand{\bcl}{\begin{center}}
\newcommand{\ecl}{\end{center}}
\newcommand{\brl}{\begin{right}}
\newcommand{\erl}{\end{right}}
\newcommand{\ben}{\begin{enumerate}}
\newcommand{\een}{\end{enumerate}}
\newcommand{\overliner}{\begin{array}}
\newcommand{\earr}{\end{array}}
\newcommand{\btab}{\begin{tabular}}
\newcommand{\etab}{\end{tabular}}
\newcommand{\bdoc}{\begin{document}}
\newcommand{\edoc}{\end{document}}
\newcommand{\beqy}{\begin{eqnarray}}
\newcommand{\eeqy}{\end{eqnarray}}

\newcommand{\beqi}{\begin{eqnarray*}}
\newcommand{\eeqi}{\end{eqnarray*}}
\newcommand{\bitem}{\begin{itemize}}

\newcommand{\eitem}{\end{itemize}}
\newcommand{\nln}{\newline}
\newcommand{\newt}{\newtheorem}
\newcommand{\pa}{\partial}
\newcommand{\re}{{I\!\!R}}
\newcommand{\ren}{\re^N}
\newcommand{\xr}{x\in\re }
\newcommand{\x}{\times}
\newcommand{\dyle}{\displaystyle}
\newcommand{\ene}{{I\!\!N}}
\newcommand{\irn}{\int\limits_{\re^N}}
\newcommand{\io}{\int\limits_{\O}}
\newcommand{\meas}{{\rm meas\,}}

\newcommand{\sign}{{\rm sign}}
\newcommand{\map}{\longrightarrow }
\newcommand{\imp}{\Longrightarrow }
\renewcommand{\div}{\nabla\cdot }
\newcommand{\sen}{{\rm sen\,}}
\newcommand{\tg}{{\rm tg\,}}
\newcommand{\arcsen}{{\rm arcsen\,}}
\newcommand{\arctg}{{\rm arctg\,}}
\newcommand{\supp}{{\textsl supp\ }}
\newcommand{\ity}{\int_{-\iy}^{+\iy}}
\newcommand{\limit}{\lim\limits}
\newcommand{\limi}{\limit_{n\to\infty}}
\newcommand{\sumi}{\sum\limits_{n=1}^{\infty}}
\newcommand{\ulu}{\underline u}
\newcommand{\ulw}{\underline w}
\newcommand{\ulz}{\underline z}
\newcommand{\ulv}{\underline v}
\newcommand{\uls}{\underline s}
\newcommand{\olu}{\overline u}
\newcommand{\olv}{\overline v}
\newcommand{\ols}{\overline s}
\newcommand{\ob}{\overline\b}
\newcommand{\ovar}{\overline\var}
\newcommand{\wv}{\widetilde v}
\newcommand{\wu}{\widetilde u}
\newcommand{\ws}{\widetilde s}
\renewcommand{\a }{\alpha }
\renewcommand{\b }{\beta }
\newcommand{\g }{\gamma}
\newcommand{\G }{\Gamma }
\renewcommand{\d }{\delta }

\newcommand{\D }{\Delta }
\newcommand{\e }{\varepsilon }
\newcommand{\z }{\zeta }
\renewcommand{\l }{\lambda }
\renewcommand{\L }{\Lambda }
\newcommand{\m }{\mu }
\newcommand{\n }{\nabla }
\newcommand{\s }{\sigma }
\newcommand{\Sig }{\Sigma }
\renewcommand{\t }{\tau }
\newcommand{\var }{\varphi }
\renewcommand{\o }{\omega }
\renewcommand{\O }{\Omega }
\newcommand{\bR}{{\bf R}}
\newcommand{\bC}{{\bf C}}
\newcommand{\bZ}{{\bf Z}}
\newcommand{\bN}{{\bf N}}
\newcommand{\bQ}{{\bf Q}}
\newcommand{\bK}{{\bf K}}
\newcommand{\bI}{{\bf I}}
\newcommand{\bv}{{\bf v}}
\newcommand{\bV}{{\bf V}}
\DeclareMathOperator{\suppo}{supp} \DeclareMathOperator{\di}{div}
\def\qed{\unskip\kern 6pt \penalty 500
\raise -2pt\hbox{\vrule \vbox to10pt{\hrule width 4pt
\vfill\hrule}\vrule}\par}

\newenvironment{Proof}{\removelastskip\vskip12pt
plus 1pt \noindent\em\rm}{\hfill {\qed \hskip .2cm}}

\title{Pointwise estimates for solutions of \\ semilinear parabolic inequalities with a potential}
\author{Luigi Montoro
\thanks{Dipartimento di Matematica e Informatica,
Universit\`a della Calabria, via Pietro Bucci, cubo 31B, 87036
Rende (CS), Italy (montoro@mat.unical.it).}\;,
Fabio
Punzo \thanks{Dipartimento di Matematica e Informatica,
Universit\`a della Calabria, via Pietro Bucci, cubo 31B, 87036
Rende (CS), Italy (fabio.punzo@unical.it).}\; and Berardino Sciunzi
\thanks{Dipartimento di Matematica e Informatica,
Universit\`a della Calabria, via Pietro Bucci, cubo 31B, 87036
Rende (CS), Italy (sciunzi@mat.unical.it).}}

\date{}

\begin{document}
\maketitle
\begin{abstract}
\noindent We obtain pointwise estimates for solutions of semilinear parabolic equations with a potential on connected domains both of $\mathbb R^n$ and of general Riemannian manifolds.
\end{abstract}

%\subjclass{35K61, 35K67, 35B99, 35B40, 35B51.}

\bigskip
\smallskip

\section{Introduction}\setcounter{equation}{0}
We are concerned with solutions of semilinear parabolic equations of the following type:
\begin{equation}\label{e1}
\pa_t u - \Delta u \, + V u^q=\,  f \quad \textrm{in}\;\; Q_T:=\Omega\times (0, T]\,,
\end{equation}
%\begin{equation}\label{e1}
%\left\{
%\begin{array}{ll}
%\,   \pa_t u -\Delta u \,=\, V u^q + f
%&\textrm{in}\,\,\Omega\times (0, T]=:Q_T
%\\& \\
%\textrm{ }u \, = u_0& \textrm{in\ \ } \Omega\times \{0\} \,,
%\end{array}
%\right.
%\end{equation}
where $\Omega\subseteq M$ is a connected domain on a complete Riemannian manifold, the potential $V=V(x,t)$ and the source term $f=f(x,t)$ are given continuous functions in $Q_T$. Moreover, we suppose that $f\geq 0$, $f\not \equiv 0$, while $V$ can be signed. We consider both the case $q>0$ and $u\geq 0$, and that $q<0$ and $u>0$.

\smallskip

The elliptic counterpart of equation \eqref{e1}, that is
\begin{equation}\label{e1i}
-\Delta u + V u^q = f \quad \textrm{in}\;\; \Omega,
\end{equation}
with $V$ and $f$ continuous functions defined in $\Omega$, has been largely investigated in the literature.
In particular, in \cite{GrigV} pointwise estimates for the solutions of \eqref{e1i} have been obtained. Indeed, in \cite{GrigV} also more general divergence form elliptic operators with smooth coefficients have been addressed.
Assume that the Dirichlet Green function of $-\Delta$ in $\Omega$ exists, and denote it by $G^\Omega(x,y).$
Set
$$H(x)\,:=\, \int_\Omega G^\Omega(x,y) f(y) d\mu(y)\,;$$
assume that $H(x)<\infty$ for all $x\in \Omega$, and that 
$$\tilde H (x)\,:=\,\int_\Omega G^\Omega(x,y) h^q(y) V(y) d\mu(y)$$
is well-defined. In \cite{GrigV} it is shown that if $q>0$, then $u$ satisfies a pointwise estimate from below, in terms of the functions $H$ and $\tilde H$. On the other hand, if $q<0$, then $u$ satisfies a similar pointwise estimate from above.  Moreover,  using similar inequalities, sufficient conditions for the existence of positive solutions of equation \eqref{e1i} have been obtained, provided $\Omega$ is relatively compact. Observe that in particular cases the results established in \cite{GrigV} have been already shown in the literature (see, e.g., \cite{BC}, \cite{BK}, \cite{FrV}, \cite{GrigH2}, \cite{GrigH}, \cite{KV})\,. However, in the remarkable paper \cite{GrigV} it is given a unified approach for treating all the values of $q\in \mathbb R\setminus\{0\}$, a general signed potential $V$, and a general divergence form operator, also on domains of Riemannian manifolds.

\smallskip

Recently, also parabolic equations with a potential on Riemannian manifolds have been investigated (see, e.g.,  \cite{BPT}, \cite{MMP}, \cite{Pu1}, \cite{Zhang}); however, it seems that in general pointwise estimates for solutions of equation \eqref{e1} have not been addressed.  In this paper we aim at obtaining pointwise estimates for solutions of \eqref{e1}, in the same spirit of the results in \cite{GrigV}, concerning elliptic equations.

Let $p$ the {\it heat kernek} in $\Omega$ (see Section \ref{mf}); for any $f\in C(Q_T)$, define for all $(x,t)\in Q_T$
\begin{equation}\label{e20}
\mathcal S^{\Omega}[f](x,t):=\int_0^t \int_\Omega p(x,y, t-s) f(y,s) d\mu(y) ds\,,
\end{equation}
provided that
\begin{equation}\label{d71}
 \int_0^t \int_\Omega p(x,y, t-s) \big|f(y,s)\big| d\mu(y) ds<\infty\quad \textrm{for every}\;\; x\in \Omega, t\in (0, T]\,.
\end{equation}
Furthermore, for any $u_0\in C(\Omega)\cap L^\infty(\Omega), u_0\geq 0$, define for all $(x,t)\in Q_T$
\begin{equation}\label{m10}
\mathcal R^{\Omega}[f; u_0](x,t):=\mathcal S^{\Omega}[f](x,t) + \int_\Omega p(x,y, t) u_0(y) d\mu(y)\,.
\end{equation}

\smallskip

We prove that for $q>0$ any solution of problem
\begin{equation}\label{e30}
\left\{
\begin{array}{ll}
\,  \pa_t u -\Delta u + V u^q \geq f\,, \;\; u \geq 0\,,
&\textrm{in}\,\, Q_T\,
\\& \\
\textrm{ } u\geq u_0 & \textrm{in}\,\, \Omega\times\{0\}
\end{array}
\right.
\end{equation}
satisfies a certain pointwise estimate from below in terms of the functions $\mathcal R^\Omega[f; u_0]$ and $\mathcal S^\Omega[h^q V]$, provided that $\mathcal S^\Omega\big[h^q |V|\big]<\infty$ in $Q_T$, where
\begin{equation}\label{e31b}
h:= \mathcal R^{\Omega}[f; u_0]\quad \textrm{in}\;\; Q_T\,.
\end{equation}
Moreover, if $q<0$, then for any solution of problem
\begin{equation}\label{e31}
\left\{
\begin{array}{ll}
\,  \pa_t u -\Delta u + V u^q \leq f\,, \;\; u > 0\,,
&\textrm{in}\,\, Q_T\,\\& \\
\textrm{ } u=0 & \textrm{in}\,\, \pa\Omega\times(0, T]
\\& \\
\textrm{ } u\leq u_0 & \textrm{in}\,\, \Omega\times\{0\}\,,
\end{array}
\right.
\end{equation}
a similar estimate from above is obtained. Indeed, note that in the case $q<0$, as well as in the elliptic case, a suitable extra pointwise condition at {\it infinity} for the solution is required. However, in the parabolic case,  if $M$ is stochastically complete, such a condition can be replaced by a growth condition at {\it infinity}, which is a weaker assumption.

Moreover, when $\Omega$ is relatively compact, we give sufficient conditions for existence of positive solutions of problem
\begin{equation}\label{d39}
\left\{
\begin{array}{ll}
\,   \pa_t u -\Delta u \,+ V u^q =\, f
&\textrm{in}\,\,Q_T
\\& \\
\textrm{ }u \, = 0& \textrm{in\ \ } \pa\Omega\times (0, T]
\\& \\
\textrm{ }u \, = u_0& \textrm{in\ \ } \Omega\times \{0\} \,,
\end{array}
\right.
\end{equation}
that are based on estimates analogous to those described above. We should note that our results seem to be new also in the case that $M=\mathbb R^n$.

\smallskip

In order to prove our results, we adapt to parabolic equations the methods used in \cite{GrigV}.  At first we prove our pointwise estimates assuming that $\Omega$ is a relatively compact connected domains, and replacing $h$ defined in \eqref{e31b} by a function
$\zeta \in C^{2,1}(Q_T)\cap C(\bar Q_T)$ that satisfies
\begin{equation}\label{e38}
\pa_t \zeta - \Delta \zeta \geq 0 \quad \textrm{in}\;\; Q_T\,,
\end{equation}
\begin{equation}\label{e37}
\zeta>0\quad \textrm{in}\;\; \Omega\times [0, T]\,.
\end{equation}
To do that the main step is to consider the equation solved by $u v$, where
\[v:= \phi^{-1}\left(\frac u h \right),\]
$\phi$ being an appropriate smooth function.
Then a suitable approximation procedure is used to obtain the desired estimates in possible not relatively compact domain $\Omega$, with $h$ defined in \eqref{e31b}.
In our arguments a special role is played by an appropriate comparison result, that is applied to the function $uv$. Note that the proof of such a comparison result is quite different from that in \cite{GrigV} for the elliptic case. Furthermore, on a special class of Riemannian manifolds, including the stochastically completes ones, we can show a refined comparison result. In view of this, we can show the estimates from above in the case $q<0$, only assuming growth conditions at infinity on the solutions of \eqref{e1}.

\smallskip

The paper is organized as follows. In Section \ref{mf} we recall some basic notions in Riemannian Geometry and in Analysis on manifolds that will be used in the sequel. Then we state our main results in Section \ref{sr}. In Section \ref{ar} we show some preliminary results, including the comparison results mentioned above, that will be essential in the proofs of the main theorems, that can be found in Sections \ref{dim} and \ref{dim2}.

\section{Mathematical framework}\label{mf}
\setcounter{equation}{0}
Let $M$ be an $n-$dimensional Riemannian manifold with a Riemannian metric tensor $g=(g_{ij})$. In any chart with coordinates $x_1,  x_2\ldots, x_n$, the associated Laplace-Beltrami operator is given by
\[ \Delta u = \frac 1{\sqrt{\textrm{det} g}}\sum_{i,j=1}^n \pa_{x_i}\big(\sqrt{\textrm{det} g} g^{ij} \pa_{x_j} u \big)\,,\]
where $\textrm{det} g$ is the determinant of the matrix $g=(g_{ij}), (g^{ij})$ is the inverse matrix of $(g_{ij}),$ and $u\in C^2(M)$.  The Riemannuan measure $d\mu$ in the same chart reads by
\[ d\mu = \sqrt{\textrm{det} g} dx_1\ldots dx_n\,;\]
furthermore, the gradient of a function $u\in C^1(M)$ is
\[ (\nabla u)^i = \sum_{j=1}^n  g^{ij} \pa_{x_j} u\,\quad (i=1,\ldots, n).\]

For any $f, g\in C^2(M)$ we have
\begin{equation}\label{e2}
\Delta (f g)\,=\, f \Delta g + 2 \langle \nabla f, \nabla g \rangle + g \Delta f\,.
\end{equation}
Moreover, for any $w\in C^2(M)$ and $\phi \in C^2(\re)$ there holds
\begin{equation}\label{e3}
\Delta[\phi(w)]\,=\, \phi'(w) \Delta w + \phi''(w) |\nabla w|^2\,.
\end{equation}

\medskip

We denote by $\partial_\infty M$ the {\it infinity point} of the one-point compactification of M (see for example [19, Sec. 5.4.3]). For any function $u:\Omega\subseteq M\to \mathbb R$ we write
\[\lim_{x\to \pa_\infty M} u(x)=0\]
to indicate that $u(x)\to 0$ as $d(x,o)\to \infty$, $o\in M$ being a fixed point; here and hereafter $d(x,y)$ denotes the geodesic distance from $x$ to $y$.
Similarly we mean equalities and inequalities involving $\liminf$ and $\limsup$.

\medskip
\smallskip

By standard results (see, e.g., \cite{Grig}) the {\it heat kernel} in $\Omega$, $p(x,y,t)$, is well-defined. For each fixed $y\in \Omega$, $p(x,y,t)$ is the smallest positive solution of equation
\begin{equation}\label{d60}
\pa_t p - \Delta p = 0\quad \textrm{in}\;\; Q_T\,,
\end{equation}
such that
\[\lim_{t\to 0^+}p(x, y, t)=\delta_y\,,\]
where $\delta_y$ is the {\it Dirac delta}  concentrated at $y$.
Moreover, $p\in C^\infty(\Omega\times \Omega\times (0, \infty)),$
\[p(x,y, t)>0 \quad \textrm{for any}\;\; x,y\in \Omega, t>0,\]
\[p(x,y,t)=p(y,x,t) \quad \textrm{for any}\;\; x,y\in \Omega, t>0\,,\]
\[p(x,y,t)=\int_\Omega p(x,z,s) p(z,y, t-s) d\mu(y) \quad \textrm{for any}\;\; t>0, 0<s<t, x,y\in \Omega\,,\]
\[ \int_\Omega p(x,y,t) d\mu(y) \leq 1 \quad \textrm{for any}\;\; x\in \Omega, t>0\,.\]
Furthermore, (see \cite[Theorem 7.16]{Grig3}) for any $u_0\in C(\Omega)\cap L^\infty(\Omega)$, the function
\[ v(x,t):= \int_M p(x,y, t)u_0(y) d\mu(y)\,,\;\;\; x\in \Omega, t>0\]
belongs to $C^\infty(\Omega\times \Omega \times (0,\infty))$,
satisfies equation \eqref{d60}, and
\[ v(x,t)\to u_0(x)\quad\; \textrm{as}\,\, t\to 0^+ \quad\textrm{locally uniformly w.r.t.}\,\, x\in \Omega\,.\]
In addition, if $\pa \Omega$ is smooth, then $v\in C(\bar Q_T)$, and
\[v=0\quad \textrm{in}\;\; \pa \Omega\times (0, T]\,.\]

\smallskip
\smallskip

As usual, we say that $f$ is locally Holder continuous in $Q_T$, if there exists $\a\in (0,1)$ such that for any compact subset $K\subset \Omega, 0<\tau\leq T$
\[ | f(x,t)- f(y, s) | \leq L \big[d(x,y)^\a + |t-s|^{\frac{\alpha}2}\big]\quad \textrm{for all}\;\; x,y\in K, t, s\in (\tau, T)\,,\]
for some $L=L_{K, \tau}>0$. We set $$C^{2,1}(Q_T):=\left\{u: Q_T\to \mathbb R\,|\, \frac{\pa^2 u}{\pa x_i\pa x_j}, \pa_t u\in C(Q_T)\,\; \textrm{for any}\;\;i,j=1,\ldots, n\, \right\}\,.$$

We have that (see, e.g. \cite{Aubin}) if \eqref{d71} holds and $f$ is locally Holder continuous in $Q_T$ and $u_0\in C(\Omega)\cap L^\infty(\Omega)$, then the function $h$ defined in \eqref{e31b} satisfies $h\in C^{2,1}(Q_T)$ and
\begin{equation}\label{d72}
\pa_t u - \Delta u\,=\, f \quad \textrm{in}\;\; Q_T\,.
\end{equation}
Moreover, if $f\in L^\infty(Q_T)$ and $u_0\in C(\Omega)\cap L^\infty(\Omega)\,$, then $h\in C(\Omega\times [0, T])$ and
\[ h = u_0\quad \textrm{in}\;\; \Omega\times \{0\}\,.\]
Finally, if $\pa \Omega$ is smooth and $f\in C(\bar Q_T)$, then
\[ h = 0 \quad \textrm{in}\;\; \pa \Omega\times (0, T]\,.\]

\section{Statements of the main results}\label{sr}
\setcounter{equation}{0}
Set
$$\chi_u(x):=\left\{
\begin{array}{ll}
\,   1
&\textrm{if}\,\,u(x)>0
\\& \\
\textrm{ }0 & \textrm{if\ \ } u(x)<0 \,.
\end{array}
\right.
$$

We can prove the pointwise estimates for solutions of \eqref{e1} contained in the following theorem.
\begin{theorem}\label{teo1}
Let $\Omega\subseteq M$ be an open connected subset. Suppose that $V, f\in C(Q_T), f\geq 0, f\not\equiv 0$ in $Q_T$, $u_0\in C(\Omega)\cap L^\infty(\Omega), u_0\geq 0$. Assume that $u\in C^{2,1}(Q_T)\cap C
(\bar Q_T)$ satisfies \eqref{e30} if $q>0$,
or that $u$ satisfies problem \eqref{e31} and
\begin{equation}\label{e49}
\lim_{x\to \pa_\infty M}\sup_{t\in (0, T]} u(x,t)=0\,,
\end{equation}
if $q<0$.
Let  \eqref{d71} be satisfied, and let $h$ be defined by \eqref{e31b}. Moreover, assume that
$$\mathcal S^{\Omega}[h^q |V|](x,t)<\infty \quad \textrm{for all}\;\; (x,t)\in Q_T,\;\; \textrm{if}\;\; q<0\,\, \textrm{or}\;\; q\geq 1,$$
or that
\begin{equation}\label{e31c}
\mathcal S^{\Omega}[\chi_u h^q |V|](x,t)<\infty \quad \textrm{for all}\;\; (x,t)\in Q_T,\;\; \textrm{if}\;\; 0<q<1\,.
\end{equation}

Then the following statements hold for all $(x,t)\in Q_T$.

\smallskip

\noindent $(i)$ If $q=1$, then
\begin{equation}\label{e32}
u(x,t) \geq h(x,t) e^{-\frac 1{h(x,t)} \mathcal S^{\Omega}[h V](x,t)}\,.
\end{equation}

\noindent $(ii)$ If $q>1$, then
\begin{equation}\label{e33}
-(q-1) \mathcal S^{\Omega}[h^q V](x,t) < h(x,t)\,,
\end{equation}
and
\begin{equation}\label{e34}
u(x,t) \geq \frac{h(x,t)}{\left\{1+(q-1)\frac{\mathcal S^{\Omega}[h^q V](x,t)}{h(x,t)} \right\}^{\frac 1{q-1}}}\,.
\end{equation}

\noindent $(iii)$ If $0<q<1$, then
\begin{equation}\label{e35}
u(x,t) \geq h(x,t)\left\{1-(q-1)\frac{\mathcal S^{\Omega}[\chi_u h^q V](x,t)}{h(x,t)}\right\}^{\frac 1{1-q}}_+\,.
\end{equation}

\noindent $(iv)$ If $q<0$, then \eqref{e33} holds, and
\begin{equation}\label{e36}
u(x,t) \leq h(x,t) \left\{1-(1-q)\frac{\mathcal S^{\Omega}[h^q V](x,t)}{h(x,t)}\right\}^{\frac 1{1-q}}\,,
\end{equation}
\end{theorem}

Furthermore, in the case that $f\equiv 0$, we can prove the following estimates.

\begin{theorem}\label{teo3}
Let $\Omega\subseteq M$ be an open connected subset. Let $V\in C(Q_T).$ Suppose that $u\in C^{2,1}(Q_T)$ satisfies either
\begin{equation}\label{e41}
\pa_t u - \Delta u + V u^q \geq 0,\;\; u\geq 0\;\, \textrm{in}\;\; Q_T,\,\,\, \textrm{if}\,\; q>0\,,
\end{equation}
or
\begin{equation}\label{e42}
\pa_t u - \Delta u + V u^q \leq 0, \;\; u>0\;\, \textrm{in}\;\,\; Q_T,\,\,\, \textrm{if}\,\; q<0\,.
\end{equation}
Moreover, assume that
\[ \mathcal S^{\Omega}[|V|](x,t)<\infty \quad \textrm{for all}\;\; (x,t)\in Q_T,\;\; \textrm{if}\;\; q<0\,\, \textrm{or}\;\; q\geq 1,
\]
or that
\begin{equation}\label{h31}
 \mathcal S^{\Omega}[\chi_u |V|](x,t)<\infty \quad \textrm{for all}\;\; (x,t)\in Q_T,\;\;\, \textrm{if}\;\; 0<q<1\,,
\end{equation}

Then the following statements hold\,.

\smallskip

\noindent $(i)$ If $q=1, u\in C(\bar Q_T),$
\begin{equation}\label{e43a}
u\geq 1\quad \textrm{in}\;\; \big[\pa \Omega\times (0, T] \big]\cup\big[\Omega\times \{0\}\big]\,,
\end{equation}
\begin{equation}\label{e43}
 \liminf_{x\to \pa_\infty M}\inf_{t\in (0, T]} u(x,t) \geq 1\,,
\end{equation}
then
\begin{equation}\label{e44}
u(x,t)\geq e^{-\mathcal S^{\Omega}[V](x,t)}\quad \textrm{for all}\;\; (x,t)\in Q_T\,.
\end{equation}

\noindent $(ii)$ If $q>1$ and
\begin{equation}\label{e45}
\lim_{t\to 0^+}\inf_{x\in \Omega} u(x,t) =\infty,
\quad \lim_{d(x,\partial\Omega)\to 0} \inf_{t\in (0, T]} u(x,t)=\infty\,,
\quad  \lim_{x\to \pa_\infty M} \inf_{t\in (0, T]} u(x,t)=\infty\,,
\end{equation}
then
\begin{equation}\label{e46}
\mathcal S^{\Omega}[V](x,t)>0\,,
\end{equation}
and
\begin{equation}\label{e47}
u(x,t)\geq \left\{(q-1)\mathcal S^{\Omega}[V](x,t)\right\}^{-\frac 1{q-1}}\,.
\end{equation}

\noindent $(iii)$ If $0<q<1$, then
\begin{equation}\label{e48}
u(x,t) \geq \left\{-(1-q)\mathcal S^{\Omega}[\chi_u V](x,t)\right\}^{\frac 1{q-1}}_+\,.
\end{equation}

\noindent $(iv)$ If $q<0, u\in C(\bar Q_T)$,
\begin{equation}\label{e49a}
u=0 \quad \textrm{in}\;\; \big[\pa\Omega\times (0, T]\big] \cup \big[ \Omega\times\{0\}\big]\,,
\end{equation}
and
\eqref{e49} is satisfied,
then
\begin{equation}\label{e50h}
\mathcal S^{\Omega}[V](x,t)<0\,,
\end{equation}
and
\begin{equation}\label{e51}
u(x,t) \leq \left\{-(1-q)\mathcal S^{\Omega}[V](x,t)\right\}^{\frac 1{q-1}}\,.
\end{equation}
\end{theorem}

\smallskip

In the next theorem, we give sufficient conditions for the existence of nonnegative solutions of problem \eqref{d39}, in the case that $\Omega$ is relatively compact, and $u_0\in C(\bar \Omega)$, with $u_0=0$ on $\partial \Omega$. 
Note that, the last compatibility condition allows us to construct solutions that attain continuously zero on the whole parabolic boundary.
Moreover, we establish two-sided pointwise estimates for such solutions.

\begin{theorem}\label{teo4}
Let $\Omega\subset M$ be a connected relatively compact subset with  boundary $\pa \Omega$ of class $C^1$. Suppose that $f$ and $V$ are locally Holder continuous in $Q_T$, and that $f\in C(\bar Q_T), f\geq 0, f\not \equiv 0$.  Assume that $u_0\in C(\bar\Omega)\,, u_0=0$ on $\partial \Omega$. Let \eqref{d71} be satisfied,  and let $h$ be defined by \eqref{e31b}. Then the following statements hold.

\smallskip

\noindent $(i)$ Suppose that $q>1, V\leq 0$, and that
\begin{equation}\label{d35}
-\mathcal S^\Omega[h^q V](x,t) \leq \left(1-\frac 1 q \right)^q \frac 1{q-1} h(x,t) \quad \textrm{for all}\;\; (x,t)\in Q_T\,.
\end{equation}
 Then a nonnegative solution $u\in C^{2,1}(Q_T)\cap C(\bar Q_T)$ of problem \eqref{d39} exists; moreover,
\begin{equation}\label{d36}
\frac{h(x,t)}{\left\{1+(q-1)\frac{\mathcal S^\Omega[h^q V](x,t)}{h(x,t)}\right\}^{\frac 1{q-1}}}\leq u(x,t) \leq \frac q{q-1}h(x,t)\quad \textrm{for all}\;\; (x,t)\in Q_T\,.
\end{equation}

\smallskip

\noindent $(ii)$ Suppose that $q<0, V\geq 0,$ and that
\begin{equation}\label{d37}
\mathcal S^\Omega[h^q V](x,t)\leq \left(1- \frac 1 q \right)^{q}\frac 1{1-q} h(x,t)\quad \textrm{for all}\;\; (x,t)\in Q_T\,.
\end{equation}
Then a positive solution $u\in C^{2,1}(Q_T)\cap C(\bar Q_T)$ of problem \eqref{d39} exists; moreover,
\begin{equation}\label{d38}
\frac 1{1-\frac 1 q} h(x,t) \leq u(x,t) \leq \left\{1-(1-q)\frac{\mathcal S^\Omega[h^q V](x,t)}{h(x,t)}\right\}^{\frac 1 {1-q}} h(x,t)\quad \textrm{for all}\;\; (x,t)\in Q_T\,.
\end{equation}
\end{theorem}

\subsection{Further results for $q<0$}
Consider domains $\Omega$ that are not relatively compact. If $q<0$, under suitable hypotheses, we can remove condition  \eqref{e49} and then getting Theorem \ref{teo1}-$(iv)$ and in Theorem \ref{teo3}-$(iv)$.

\medskip

We assume that there exist $\mu>0$ and a subsolution $Z$ of equation
\begin{equation}\label{h1}
\Delta Z \,=\, \mu Z \quad \textrm{in}\;\; \Omega\,,
\end{equation}
such that
\begin{equation}\label{h3}
\sup_\Omega Z <\infty\,,\quad \lim_{x\to \pa_\infty M} Z(x)\,=\, -\infty\,.
\end{equation}
By a subsolution of \eqref{h1} we mean a function $Z\in C^2(\Omega)$ such that 
\begin{equation}\label{asildhfg}
\Delta Z \,\geq\, \mu Z\quad \textrm{in}\;\; \Omega\,.
\end{equation}
Observe that our results remain true if $Z$ is continuous in $\Omega$ and satisfies \eqref{asildhfg} in the distributional sense.
Note that, in the case $\Omega=M$, the existence of such a subsolution $Z$ implies that $M$ is stochastically complete (see \cite{Grig}),
i.e. 
\begin{equation}\nonumber
\int_M p(x,y,t)\,d\mu(y)\,=1\qquad\text{for all}\quad x\in M,\quad t>0\,.
\end{equation}
We refer the reader to \cite{Grig} for sufficient and necessary condition for the existence of such subsolution $Z$. We limit ourselves to observe that
such a subsolution $Z$ exists 
for instance on $\mathbb{R}^n$, $n\geq 3$, and on the hyperbolic space $\mathbb{H}^n$, $n\geq 2$.
\begin{theorem}\label{teo5}
Let $q<0.$ Let $\Omega\subseteq M$ be an open not relatively compact connected subset. Suppose that $V, f\in C(Q_T), f\geq 0, f\not\equiv 0$ in $Q_T$, $u_0\in C(\Omega)\cap L^\infty(\Omega), u_0\geq 0$. Assume that $u\in C^{2,1}(Q_T)\cap C
(\bar Q_T)$ satisfies \eqref{e31}. Let conditions \eqref{d71} and \eqref{e31c} be satisfied, and let $h$ be defined by \eqref{e31b}.
Let there exist $\mu>0$ and a subsolution $Z$ of equation \eqref{h1}, which satisfies \eqref{h3}.
Moreover, suppose that
\begin{equation}\label{h30}
\limsup_{x\to \pa_\infty M} \frac{\sup_{t\in (0, T]}h^q(x,t)[u^{1-q}(x,t) - h^{1-q}(x,t)]}{|Z(x)|}\,\leq\,0\,.
\end{equation}
Then \eqref{e33} and \eqref{e36} hold.
\end{theorem}

\begin{theorem}\label{teo6}
Let $q<0$. Let $\Omega\subseteq M$ be an open not relatively compact connected subset. Let $V\in C(Q_T).$ Suppose that $u\in C^{2,1}(Q_T)\cap C(\bar Q_T)$ satisfies \eqref{e42} and \eqref{e49a}. Let condition \eqref{h31} be satisfied. Let there exist $\mu>0$ and a subsolution $Z$ of equation \eqref{h1}, which satisfies \eqref{h3}. Moreover, suppose that
\begin{equation}\label{h32}
\limsup_{x\to \pa_\infty M} \frac{\sup_{t\in (0, T]}u^{1-q}(x,t)}{|Z(x)|}\,\leq\,0\,.
\end{equation}
Then \eqref{e50h} and \eqref{e51} hold.
\end{theorem}

\begin{remark}
It is easily seen that both condition \eqref{h30} and \eqref{h32} are weaker than condition \eqref{h49}.
\end{remark}

\section{Auxiliary results}\label{ar}\setcounter{equation}{0}
This section is devoted to some preliminary results that will be used to prove Theorems \ref{teo1}, \ref{teo3}, \ref{teo4}\,.

\begin{lemma}\label{l1}
Let $v, h\in C^{2,1}(Q_T), \phi\in C^2(I)$ with $v(Q_T)\subseteq I, I$ being an interval in $\re$. Then
\begin{equation}\label{e4}
\begin{aligned}
& \pa_t[h \phi(v)] -\Delta[h \phi(v)] \\ &\,=\,  \phi'(v)[\pa_t(h v) - \Delta (h v)] - \phi''(v) |\nabla v|^2 h + [\phi(v) - v \phi'(v)](\pa_t h - \Delta h)\quad \textrm{in}\;\; Q_T\,.
\end{aligned}
\end{equation}
In particular, if $\phi'\neq 0$ in $I$, then
\begin{equation}\label{e5}
\begin{aligned}
& \pa_t(h v) -\Delta(h v) \\ &\,=\,  \frac{\pa_t[h\phi(v)]- \Delta[h \phi(v)]}{\phi'(v)} + \frac{\phi''(v)}{\phi'(v)} |\nabla v|^2 h + \left(v - \frac{\phi(v)}{\phi'(v)}\right)(\pa_t h - \Delta h)\quad \textrm{in}\;\; Q_T\,.
\end{aligned}
\end{equation}
\end{lemma}

\noindent{\it Proof\,.}
Clearly,
\begin{equation}\label{e16}
\pa_t[h \phi(v)]\,=\,\phi'(v) \pa_t(hv) + [\phi(v)- v \phi'(v)] \pa_t h\,.
\end{equation}
Moreover, in view of \eqref{e2} with $f=h, g=\phi(v)$, and  in view of  \eqref{e3} with $w=v$ we get
\[\Delta[h \phi(v)]=\phi(v) \Delta h + h [\phi'(v)\Delta v + \phi''(v) |\nabla v|^2  ] + 2 \phi'(v) \langle \nabla h, \nabla v\rangle\,. \]
Thus
\begin{equation}\label{e17}
\begin{aligned}
 \Delta[h \phi(v)]  =\, &  \phi'(v) \Delta(h v) \\
& + \phi''(v) |\nabla v|^2h + [\phi(v)- v \phi'(v)] \Delta h\,.
\end{aligned}
\end{equation}
>From \eqref{e16} and \eqref{e17} we easily obtain \eqref{e4}, and then \eqref{e5}.

\hfill $\square$

\begin{lemma}\label{l2}
Let $I\subseteq \re$ be an interval. Let $\phi\in C^2(I), \phi>0, \phi'>0$ in $I$. Let $v, h\in C^{2,1}(Q_T) $ with $h>0, v(\Omega)\subseteq I$. Set
\[u:= h \phi(v)\,.\]
Let $V\in C(Q_T), q\in \re\setminus\{0\}$. If
\begin{equation}\label{e6}
\pa_t u - \Delta u + V u^q \geq \pa_t h - \Delta h \quad \textrm{in}\;\; Q_T\,,
\end{equation}
then
\begin{equation}\label{e7}
\begin{aligned}
& \pa_t(h v)\,-\, \Delta (hv)\, + h^q V \frac{\phi(v)^q}{\phi'(v)}\, \\
& \, \geq\, \left(v- \frac{\phi(v)-1}{\phi'(v)}\right)(\pa_t h- \Delta h ) + \frac{\phi''(v)}{\phi'(v)} |\nabla v|^2 h\quad \textrm{in}\;\, Q_T\,.
\end{aligned}
\end{equation}
If
\begin{equation}\label{e8}
\pa_t u - \Delta u + V u^q \leq \pa_t h - \Delta h \quad \textrm{in}\;\; Q_T\,,
\end{equation}
then
\begin{equation}\label{e9}
\begin{aligned}
& \pa_t(h v)\,-\, \Delta (hv)\, + h^q V \frac{\phi^q(v)}{\phi'(v)}\, \\
& \, \leq\, \left(v- \frac{\phi(v)-1}{\phi'(v)}\right)(\pa_t h- \Delta h ) + \frac{\phi''(v)}{\phi'(v)} |\nabla v|^2 h\quad \textrm{in}\;\, Q_T\,.
\end{aligned}
\end{equation}
\end{lemma}

\noindent{\it Proof\,.} From \eqref{e6} with $u= h\phi(v)$ it follows that
\begin{equation}\label{e18}
\pa_t[ h\phi(v)] - \Delta[h \phi(v)] \geq - V h^q \phi(v)^q + \pa_t h - \Delta h \,.
\end{equation}
Therefore, by \eqref{e5} and \eqref{e18},
\[
\begin{aligned}
\pa_t( h v) - \Delta(h v) & \\ & \geq - V h^q \frac{\phi(v)^q}{\phi'(v)} + \frac{\phi''(v)}{\phi'(v)} |\nabla v|^2 h+ \frac{1+  v \phi'(v) - \phi(v)}{\phi'(v)} (\pa_t h - \Delta h)\,.
\end{aligned}
\]
So, \eqref{e7} follows. The second claim can be proved in the same way. \hfill $\square$

\begin{lemma}\label{l3}
Let assumptions of Lemma \ref{l2} be satisfied. Moreover, suppose that $0\in I$, and that
\begin{equation}\label{e10}
\pa_t h - \Delta h \geq 0\quad \textrm{in}\;\;\, Q_T\,.
\end{equation}
If
\begin{equation}\label{e11}
\phi(0)=1,
\end{equation}
\begin{equation}\label{e12}
\phi'>0,\;\, \phi''\geq 0\quad \textrm{in}\;\; I\,,
\end{equation}
then
\begin{equation}\label{e13}
\pa_t(hv) - \Delta (hv) + h^q V \frac{\phi(v)^q}{\phi'(v)} \geq 0 \quad \textrm{in}\,\; Q_T\,.
\end{equation}
If \eqref{e11} holds, and
\begin{equation}\label{e14}
\phi'>0,\;\, \phi''\leq 0\quad \textrm{in}\;\; I\,,
\end{equation}
then
\begin{equation}\label{e15}
\pa_t(hv) - \Delta (hv) + h^q V \frac{\phi(v)^q}{\phi'(v)} \leq 0 \quad \textrm{in}\,\; Q_T\,.
\end{equation}
\end{lemma}

\noindent{\it Proof\,.} It is direct to see that \eqref{e11} and \eqref{e12} imply that
\begin{equation}\label{e19}
v- \frac{\phi(v)-1}{\phi'(v)} \geq 0 \quad \textrm{for all}\;\, v\in I\,.
\end{equation}
>From \eqref{e7}, \eqref{e10} and \eqref{e19} we obtain \eqref{e13}. Inequality \eqref{e15} can be deduced similarly. \hfill $\square$

\begin{remark}\label{ossl3}
Note that if $\pa_t h - \Delta h =0$ in $Q_T$, then in Lemma \ref{l3} condition \eqref{e11} can be removed.
\end{remark}

\medskip
\medskip

In the sequel, we often use the next comparison result.
\begin{proposition}\label{prop1}
Let $\Omega\subset M$ be an open subset. Assume that $g \in C(Q_T)$, and that
\begin{equation}\label{e64}
\mathcal S[|g|]<\infty\,\quad \textrm{in}\;\; Q_T\,.
\end{equation}
Let $v\in C^2(Q_T)\cap C(\bar Q_T)$ be a supersolution of problem
\begin{equation}\label{e21}
\left\{
\begin{array}{ll}
\,   \pa_t v -\Delta v \,= \, g
&\textrm{in}\,\,Q_T
\\& \\
\textrm{ } v \, = 0 & \textrm{in\ \ } \pa \Omega\times (0, T] \\&\\
\textrm{ }v \, = 0& \textrm{in\ \ } \Omega\times \{0\} \,.
\end{array}
\right.
\end{equation}
Furthermore, if $\Omega$ is not relatively compact, suppose that
\begin{equation}\label{e22}
\liminf_{x\to \pa_\infty M} \inf_{t\in (0, T]} v(x,t) \geq 0\,.
\end{equation}
Then
\begin{equation}\label{e23}
v(x,t) \geq \mathcal S^{\Omega}[g](x, t)\quad \textrm{for every}\,\, x\in \Omega, t\in [0, T]\,.
\end{equation}
\end{proposition}

\noindent{\it Proof\,.} Choose a sequence of functions $\{g_n\}$ such that $g_n$ is locally Lipschitz continuous in $Q_T$ for every $n\in \mathbb N$,
\begin{equation}\label{e50}
g_n \leq g\,,\;\; g_n\leq g_{n+1}\quad \textrm{in}\;\; Q_T\quad \textrm{for every}\,\, n\in \mathbb N\,;
\end{equation}
\begin{equation}\label{e50b}
g_n\to g \quad \textrm{in}\;\; Q_T\,\, \textrm{as}\,\, n\to \infty\,.
\end{equation}
Let us only consider the case when $\Omega$ is not relatively compact; the case when $\Omega$ is relatively compact is easier and it will be omitted.\\
\noindent Let $k\in\mathbb{N}$ that will be taken arbitrary large later on.
Fixed a point $o\in M$, by \eqref{e22}, we find a radius $R_k$ such that
\begin{equation}\label{e91}
v\geq -\frac{1}{k}\qquad \text{on}\quad  \left(\Omega \cap \partial B_{R_k} (o)\right)\times (0, T]\,.
\end{equation}
Since $v\in  C(\bar Q_T)$ we can therefore take $\Omega_k\subseteq \Omega \cap  B_{R_k} (o)$ so that
\begin{equation}\nonumber
v\geq -\frac{1}{k}\qquad \text{on}\quad \partial\Omega_k\times (0, T]\,.
\end{equation}
For each $k$ fixed, the construction of $\Omega_k$ can be carried out just observing that $v$ is uniformly continuous in $\Omega \cap B_{R_k} (o)$ and exploiting the boundary datum.
With no loss of generality we may and do assume that $R_k\rightarrow\infty$, $\Omega_k$ is smooth and
\begin{equation}\label{kdjfvhkdds}
\underset{k\in\mathbb{N}}{\cup}\,\Omega_k\,=\,\Omega\,.
\end{equation}
Therefore, by construction, we have that
 $v$ is a supersolution of the problem
\begin{equation}\label{e27}
\left\{
\begin{array}{ll}
\,   \pa_t v -\Delta v \,= \, g_n
&\textrm{in}\,\, \Omega_k\times (0, T]
\\& \\
\textrm{ } v \, \geq -k^{-1} & \textrm{in\ \ } \pa \Omega_k\times (0, T] \\&\\
\textrm{ }v \, \geq -k^{-1} & \textrm{in\ \ } \Omega_k\times \{0\} \,.
\end{array}
\right.
\end{equation}
Let now $v_{n,k}$ be the solution of the problem
\begin{equation}\label{e25}
\left\{
\begin{array}{ll}
\,   \pa_t v -\Delta v \,= \, g_n
&\textrm{in}\,\, \Omega_k\times (0, T]
\\& \\
\textrm{ } v \, = 0 & \textrm{in\ \ } \pa \Omega_k\times (0, T] \\&\\
\textrm{ }v \, = 0& \textrm{in\ \ } \Omega_k\times \{0\} \,.
\end{array}
\right.
\end{equation}
We have that
\begin{equation}\label{e53}
v_{n, k}(x,t)=\int_0^t  \int_{\Omega_k} p_k(x,y, t-s) g_n(y,s) dt d\mu(y),\quad x\in \bar \Omega_k, t\in [0, T]\,,
\end{equation}
where $p_k$ is the heat kernel in $\Omega_k$, completed with zero homogeneous Dirichlet boundary conditions.
It is known that (see, e.g., \cite{Grig}), by \eqref{kdjfvhkdds}, it follows that
\begin{equation}\label{d73}
\lim_{k\to \infty} p_k = p \quad \textrm{in}\;\; M\times M\times (0, \infty)\,.
\end{equation}
Therefore, using \eqref{e64}, \eqref{e50b} and \eqref{d73},  we can infer that
\begin{equation}\label{e26}
\lim_{n\to \infty, k\to \infty} v_{n, k} \,=\, \mathcal S^{\Omega}[g]\quad \textrm{in}\;\;  Q_T\,.
\end{equation}
On the other hand, the function $v_{n, k}-k^{-1}$ is a subsolution of problem
\begin{equation}\label{e27jjgjgjghh}
\left\{
\begin{array}{ll}
\,   \pa_t v -\Delta v \,= \, g_n
&\textrm{in}\,\, \Omega_k\times (0, T]
\\& \\
\textrm{ } v \, \leq -k^{-1} & \textrm{in\ \ } \pa \Omega_k\times (0, T] \\&\\
\textrm{ }v \, \leq -k^{-1} & \textrm{in\ \ } \Omega_k\times \{0\} \,.
\end{array}
\right.
\end{equation}
By the comparison principle, taking into account \eqref{e27} and \eqref{e27jjgjgjghh}, we deduce that
\begin{equation}\label{e28}
 v\geq v_{n, k} - k^{-1} \quad \textrm{in}\;\; \Omega_k\times [0, T]\,.
 \end{equation}
In view of \eqref{e26}, letting $k\to \infty, n\to \infty$, we obtain \eqref{e23}\,.

\hfill $\square$

We also use the next comparison result.
\begin{proposition}\label{prop1a}
Let $\Omega\subset M$ be an open subset. Assume that $g \in C(Q_T)$ and that \eqref{e64} is satisfied. Let $v\in C^2(Q_T)\cap C(\bar Q_T)$ be a subsolution of problem \eqref{e21}.
%\begin{equation}\label{e21a}
%\left\{
%\begin{array}{ll}
%\,   \pa_t v -\Delta v \,= \, g
%&\textrm{in}\,\,Q_T
%\\& \\
%\textrm{ } v \, = 0 & \textrm{in\ \ } \pa \Omega\times (0, T] \\&\\
%\textrm{ }v \, = 0& \textrm{in\ \ } \Omega\times \{0\} \,.
%\end{array}
%\right.
%\end{equation}
Furthermore, if $\Omega$ is not relatively compact, suppose that
\begin{equation}\label{e22a}
\limsup_{x\to \pa_\infty M} \sup_{t\in (0, T]} v(x,t) \leq 0\,.
\end{equation}
Then
\begin{equation}\label{e23a}
v(x,t) \leq \mathcal S^{\Omega}[g](x, t)\quad \textrm{for every}\,\, x\in \Omega, t\in [0, T]\,.
\end{equation}
\end{proposition}

The proof of Proposition \ref{prop1a} is analogous to that of Proposition \ref{prop1}; the only difference is that the sequence $\{ g_n\}$ satisfies
\begin{equation}\label{e50a}
g_n \geq g\,,\;\; g_n\geq g_{n+1}\quad \textrm{in}\;\; Q_T\quad \textrm{for every}\,\, n\in \mathbb N\,,
\end{equation}
instead of \eqref{e50}.

\medskip
\medskip

Moreover, we use the next refined comparison principles.
\begin{proposition}\label{prop1r}
Let $\Omega\subset M$ be an open, not relatively compact subset. Assume that $g \in C(Q_T)$, and that \eqref{e64} is satisfied.
Let $v\in C^2(Q_T)\cap C(\bar Q_T)$ be a subsolution of problem \eqref{e21}. Assume that there exists a  subsolution $Z$ of equation \eqref{h1} such that \eqref{h3} is satisfied.
Furthermore, suppose that
\begin{equation}\label{h10}
\limsup_{x\to \pa_\infty M}  \frac{\sup_{t\in (0, T]}v(x,t)}{|Z(x)|}\leq 0\,.
\end{equation}
Then \eqref{e23a} holds.

\end{proposition}

\noindent{\it Proof\,.}
First of all we observe that
we can assume that, for some $H>0$,
\begin{equation}\label{h2}
Z \leq - H <0 \quad \textrm{in}\;\; \Omega.
\end{equation}
In fact, if $\sup_{\Omega} Z\geq 0$, then instead of $Z$ we can consider the function
\[\tilde Z:= Z - \sup_{\Omega} Z -1, \]
that clearly satisfies \eqref{h1}, \eqref{h3} and \eqref{h2}.\\

Choose now a sequence of functions $\{g_n\}$ such that $g_n$ is locally Lipschitz continuous in $Q_T$ for every $n\in \mathbb N$,
\eqref{e50a} and \eqref{e50b} hold.
 Let $k\in\mathbb{N}$ that will be taken arbitrary large later on and fix
a point $o\in M$.
We set
\[ V_k(x,t):= -k^{-1} Z(x) e^{\mu t} \quad \big((x,t)\in Q_T\big)\,. \]
In view of \eqref{h2}, since $\mu>0$, we have that
\begin{equation}\label{h12}
V_k \geq \frac{H}{k}>0 \quad \textrm{in}\;\; Q_T\,.
\end{equation}
By \eqref{h10}, we find a radius $R_k$ such that

\begin{equation}\label{h13}
v \leq V_k \quad \textrm{in}\;\; \big(\pa B_{R_k}(o)\cap \Omega \big)\times (0, T]\,.
\end{equation}
Since $v\in  C(\bar Q_T)$ we can therefore take $\Omega_k\subseteq \Omega \cap  B_{R_k} (o)$ so that
\begin{equation}\label{fghdgjhdghjkfg}
v\leq V_k\qquad \text{on}\quad \partial\Omega_k\times (0, T]\,.
\end{equation}
With no loss of generality we may and do assume that $R_k\rightarrow\infty$, $\Omega_k$ is smooth and
\begin{equation}\label{kdjfvhkddsvvv}
\underset{k\in\mathbb{N}}{\cup}\,\Omega_k\,=\,\Omega\,.
\end{equation}
With such a construction we let $v_{n,k}$ and $p_k$
as in \eqref{e53}. It is now easy to verify that
 $V_k$ is a supersolution of the problem
\begin{equation}\label{h11}
\left\{
\begin{array}{ll}
\,   \pa_t u -\Delta u \,= \, 0
&\textrm{in}\,\,\Omega_k\times (0, T]
\\& \\
\textrm{ }  u\, = V_k & \textrm{in\ \ } \pa \Omega_k\times (0, T] \\&\\
\textrm{ }u \, =  V_k& \textrm{in\ \ } \Omega_k\times \{0\} \,.
\end{array}
\right.
\end{equation}
Inequalities \eqref{h12} and \eqref{h13} and \eqref{fghdgjhdghjkfg} easily yield that
\begin{equation}\label{h14}
v - v_{n, k}\leq V_k\quad \textrm{in}\;\; \big[ \pa \Omega_k \times (0, T]\big] \cup \big[\Omega_k\times \{0\} \big]\,.
\end{equation}

%Moreover, let $\{\zeta_R\}\subset C^\infty_c(B_R(o))$
%be a family of functions such that, for each $R>0$,
%$0\leq \zeta_R\leq 1$, $\zeta_R\equiv 1$ in $B_{R/2}(o)\,.$

Exploiting \eqref{h14} and \eqref{e50a} we can infer that  $v-v_{n, k}$ is a subsolution of problem \eqref{h11} and,
by the comparison principle, we  obtain that
\begin{equation}\label{h15}
v- v_{n, k}\leq V_k \quad \textrm{in}\;\; \Omega_k\times (0, T]\,.
\end{equation}
Letting $n\to \infty, k\to \infty$ in \eqref{h15} we deduce that
\[v\leq\mathcal S^{\Omega}[g]  \quad \textrm{in}\;\; Q_T\,.\]
 \hfill $\square$

\medskip
\medskip

Similarly, the next refined comparison principle can also be shown.
\begin{proposition}\label{prop1ar}
Let $\Omega\subset M$ be an open, not relatively compact subset. Assume that $g \in C(Q_T)$ and that \eqref{e64} is satisfied. Let $v\in C^2(Q_T)\cap C(\bar Q_T)$ be a supersolution of problem \eqref{e21}. Let there exist a subsolution $Z$ of equation \eqref{h1} such that \eqref{h3} is satisfied. Furthermore, suppose that
\begin{equation}\label{eh16}
\liminf_{x\to \pa_\infty M}\frac{ \inf_{t\in (0, T]} v(x,t)}{|Z(x)|} \geq 0\,.
\end{equation}
Then \eqref{e23} holds.
\end{proposition}

\subsection{Pointwise estimates in relatively compact domains with general smooth supersolutions}
Let $h\in C^{2,1}(Q_T)\cap C(\bar Q_T)$ be a function that satisfies \eqref{e38}, \eqref{e37}\,. Consider the following inital-boundary value inequalities
\begin{equation}\label{e39}
\left\{
\begin{array}{ll}
\,   \pa_t u -\Delta u \,+ V u^q\geq \, \pa_t h - \Delta h
&\textrm{in}\,\,Q_T
\\& \\
\textrm{ }u \, \geq  h & \textrm{in\ \ } \pa\Omega\times (0, T]
\\& \\
\textrm{ }u\geq h& \textrm{in\ \ } \Omega\times \{0\}
\\& \\
\textrm{ }u \, \geq 0& \textrm{in\ \ } Q_T\,,
\end{array}
\right. \quad (q>0)
\end{equation}
and
\begin{equation}\label{e40}
\left\{
\begin{array}{ll}
\,   \pa_t u -\Delta u \,+ V u^q\leq \, \pa_t h - \Delta h
&\textrm{in}\,\,Q_T
\\& \\
\textrm{ }u \, \leq  h & \textrm{in\ \ } \pa\Omega\times (0, T]
\\& \\
\textrm{ }u\, \leq h& \textrm{in\ \ } \Omega\times \{0\}
\\& \\
\textrm{ }u \, > 0& \textrm{in\ \ } Q_T\,.
\end{array}
\right. \quad (q<0)
\end{equation}

\medskip
\smallskip

The next result has a crucial role in the proof of Theorem \ref{teo1}. In fact, it gives the estimates \eqref{e32}-\eqref{e36}, under the extra assumption that $\Omega$ is relatively compact; moreover, a general smooth function $h$ that satisfies \eqref{e38}-\eqref{e37} is used.

\begin{theorem}\label{teo2}
Let $\Omega\subseteq M$ be a relatively compact connected subset. Let $h$ be any function belonging to $C^{2,1}(Q_T)\cap C(\bar Q_T)$ that satisfies \eqref{e38}-\eqref{e37}. Let $u\in C^{2,1}(Q_T)\cap C(\bar Q_T)$ be a solution of either \eqref{e39} or \eqref{e40}.

Moreover, assume that
$$\mathcal S^{\Omega}[h^q |V|](x,t)<\infty \quad \textrm{for all}\;\; (x,t)\in Q_T,\;\; \textrm{if}\;\; q<0\,\, \textrm{or}\;\; q\geq 1,$$ or that
$$\mathcal S^{\Omega}[\chi_u h^q |V|]<\infty \quad \textrm{for all}\;\; (x,t)\in Q_T,\;\; \textrm{if}\;\; 0<q<1\,.$$
Then \eqref{e32}--\eqref{e36} hold for all $(x,t)\in Q_T$.
\end{theorem}

\noindent{\it Proof of Theorem \ref{teo2}\,.} To begin with, we further assume  that
\begin{equation}\label{e54a}
h>0,\; u>0\;\, \textrm{in}\;\; \bar Q_T,\,\; \textrm{and}\,\, V\in C(\bar Q_T)\,.
\end{equation}
Following the proof of \cite[Theorem 3.2]{GrigV},  we choose a function $\phi$ to solve the initial value problem
\begin{equation}\label{e54}
\phi'(s)\,=\,\phi(s)^q,\;\;\, \phi(0)=1\,.
\end{equation}
For $q=1$ we have
\begin{equation}\label{e55}
\phi(s)=e^s\,,\quad s\in \mathbb R\,,
\end{equation}
while for $q\not =1$ we obtain
\begin{equation}\label{e56}
\phi(s)=[(1-q)s + 1]^{\frac 1{1-q}}\,, s\in I_q,
\end{equation}
where the interval $I_q$ is given by
\begin{equation}\label{e57}
I_q =\left\{
\begin{array}{ll}
\,  \left(-\infty, \frac 1{q-1}\right) &\textrm{if}\,\, q>1\,,
\\& \\
\textrm{ } \mathbb R & \textrm{if\ \ }\,\, q=1\,, \\&\\
\textrm{ }\left( -\frac 1{q-1}, \infty\right) & \textrm{if\ \ } q<1\,.
\end{array}
\right.
\end{equation}
There holds
\begin{equation}\label{e58}
\phi'(s)=[(1-q)s +1]^{\frac q{1-q}}\,,\;\; \phi''(s)=q[(1-q)s+1]^{\frac{2q-1}{1-q}}\,.
\end{equation}
In particular, we have
\begin{equation}\label{e59}
\phi'>0 \quad \textrm{in}\;\; I_q\,;
\end{equation}
consequently, the inverse function $\phi^{-1}:(0, \infty)\to \mathbb R$ is well-defined\,.
Moreover,
\begin{equation}\label{e60}
\phi''(s)>0 \quad \textrm{in}\;\; I_q\;\,\, \textrm{if}\,\, q>0\,,
\end{equation}
whereas
\begin{equation}\label{e61}
\phi''(s)<0 \quad \textrm{in}\;\; I_q\;\,\, \textrm{if}\,\, q<0\,.
\end{equation}
Indeed, for $0<q<1$, we extend the domain of $\phi$ to all $s\leq -\frac 1{1-q}$, by putting $\phi(s)=0,$  so that
\begin{equation}\label{e62}
\phi(s)=[(1-q)s +1]_+^{\frac 1{1-q}}\, \;\, \textrm{for all}\,\,\, s\in \mathbb R\,.
\end{equation}

Due to \eqref{e54a}, we can define
\begin{equation}\label{e63}
v:=\phi^{-1}\left(\frac{u}{h}\right)\quad \textrm{in}\,\, \bar Q_T\,;
\end{equation}
we have that $v\in C^{2,1}(Q_T)\cap C(\bar Q_T)$\,.
Let $q>0$. From \eqref{e38} and \eqref{e39} we have that the function $u= h \phi(v)$ satisfies
\begin{equation}\label{e65ab}
\pa_t u - \Delta u + V u^q\geq  \pa_t h - \Delta h \geq 0 \quad \textrm{in}\;\; Q_T\,.
\end{equation}
Thanks to \eqref{e65ab}, Lemma \ref{l3} and \eqref{e54} we get
\begin{equation}\label{e65}
\pa_t( hv) - \Delta ( hv) \geq - h^q V \quad \textrm{in}\;\; Q_T\,.
\end{equation}
Since $u\geq h$ in $\big[\pa \Omega\times (0, T] \big] \cup \big[ \Omega\times \{0\}\big],$ we have that
\begin{equation}\label{e66}
 hv = h \phi^{-1}\left(\frac u h \right)\geq h \phi^{-1}(1)=0 \quad \textrm{in}\;\, \big[\pa \Omega\times (0, T] \big] \cup \big[ \Omega\times \{0\}\big] \,.
\end{equation}
So, $h v$ is a supersolution of problem \eqref{e21} with $g=- h^q V$.  Since $\Omega$ is relatively compact, by Proposition \ref{prop1},
\begin{equation}\label{e67}
 hv \geq - \mathcal S^\Omega[h^q V]\quad \textrm{in}\;\; Q_T\,.
\end{equation}
Thus,
\begin{equation}\label{e68}
v\geq - \frac 1 h \mathcal S[h^q V]\quad \textrm{in}\;\; Q_T\,.
\end{equation}

As a consequence of \eqref{e63} and \eqref{e68} we obtain that, for $q>1$,
\begin{equation}\label{e69}
v < \frac 1{q-1},\;\;\, - h^{-1}\mathcal S^\Omega[h^q V] < \frac 1{q-1}\,.
\end{equation}
Hence, for each $q>0$, we can apply $\phi$ to both sides of \eqref{e68} to obtain
\begin{equation}\label{e70}
\frac u h \leq \phi\left(-\frac 1 h \mathcal S^\Omega[h^q V]\right)\quad \textrm{in}\;\; Q_T\,,
\end{equation}
which implies \eqref{e32}, \eqref{e34}, \eqref{e35}\,. Moreover, from \eqref{e69} it follows \eqref{e33}\,.

\smallskip

Now, assume that $q<0$. Then we have
\[\pa_t u - \Delta u + V u^q \leq \pa_t h - \Delta h \quad \textrm{in}\;\; Q_T\,.
\]
Thanks to Lemma \ref{l3} and \eqref{e15} we have
\begin{equation}\label{e65a}
\pa_t( hv) - \Delta ( hv) \leq - h^q V \quad \textrm{in}\;\; Q_T\,.
\end{equation}
Since $u\leq h$ in $\big[\pa \Omega\times (0, T] \big] \cup \big[ \Omega\times \{0\}\big],$ we have that
\begin{equation}\label{e66a}
 hv = h \phi^{-1}\left(\frac u h \right)\leq h \phi^{-1}(1)=0 \quad \textrm{in}\;\, \big[\pa \Omega\times (0, T] \big] \cup \big[ \Omega\times \{0\}\big] \,.
\end{equation}
So, $h v$ is a subsolution of problem \eqref{e21} with $g=- h^q V$.  Since $\Omega$ is bounded, by Proposition~\ref{prop1a},
\[
 hv \leq - \mathcal S^\Omega[h^q V]\quad \textrm{in}\;\; Q_T\,.
\]
Thus,
\begin{equation}\label{e68a}
v\leq - \frac 1 h \mathcal S^\Omega[h^q V]\quad \textrm{in}\;\; Q_T\,.
\end{equation}
In view of \eqref{e68a}, it follows \eqref{e33}. Moreover, applying $\phi$ to both sides of \eqref{e68a} we get
\begin{equation}\label{e70a}
\frac u h \geq \phi\left(-\frac 1 h \mathcal S^\Omega[h^q V]\right)\quad \textrm{in}\;\; Q_T\,,
\end{equation}
and then \eqref{e36}\,.

\medskip

Now we can remove the extra assumptions in \eqref{e54a}. We extend the domain $I_q$ of $\phi$ to the endpoints of $I_q$ by taking the limits of $\phi$ at the endpoints. So, the extended domain of $\phi$ is the interval
$$
\bar I_q =\left\{
\begin{array}{ll}
\,  \left[-\infty, \frac 1{q-1}\right] &\textrm{if}\,\, q>1\,,
\\& \\
\textrm{ } [-\infty, \infty] & \textrm{if\ \ }\,\, q=1\,, \\&\\
\textrm{ }\left[-\frac 1{q-1}, \infty\right] & \textrm{if\ \ } q<1\,.
\end{array}
\right.
$$
Moreover, when $0<q<1$, we extend $\phi$ to all $s\in [-\infty, \infty]$ by using \eqref{e62}\,. Hence \eqref{e32}, \eqref{e34} and \eqref{e35} can be written in the form \eqref{e70}, while \eqref{e51} in the form \eqref{e70a}.

\smallskip

Take $q>0$. Let us show \eqref{e70}. To this purpose, for every $\e>0$ set
\[ u_\e:= u+ \e\]
and define
\[ v_\e:= \phi^{-1}\left( \frac{u_\e}{h} \right) \quad \textrm{in} \;\; Q_T\,.\]
Note that since $u_\e>0$ and $h>0$ in $Q_T$, the function $v_\e$ is well-defined in $Q_T$ and $v_\e\in C^{2,1}(Q_T)$; moreover, $v_\e(Q_T)\subset I_q\,.$ From \eqref{e5} it follows that
\begin{equation}\label{e71}
\begin{aligned}
& \pa_t(h v_\e) -\Delta(h v_\e) \\ &\,=\,  \frac{\pa_t[h\phi(v_\e)]- \Delta[h \phi(v_\e)]}{\phi'(v_\e)} + \frac{\phi''(v_\e)}{\phi'(v_\e)} |\nabla v_\e|^2 h + \left(v_\e - \frac{\phi(v_\e)}{\phi'(v_\e)}\right)(\pa_t h - \Delta h)\quad \textrm{in}\;\; Q_T\,.
\end{aligned}
\end{equation}
Since
\[ \pa_t[ h \phi(v_\e)] - \Delta[h\phi(v_\e)] = \pa_t u_\e - \Delta u_\e = \pa_t u - \Delta u \quad \textrm{in}\;\; Q_T\,,  \]
we get
\begin{equation}\label{e72}
\begin{aligned}
& \pa_t(h v_\e) -\Delta(h v_\e) \\ &\,=\,  \frac{\pa_t u - \Delta u}{\phi'(v_\e)} + \frac{\phi''(v_\e)}{\phi'(v_\e)} |\nabla v_\e|^2 h + \left(v_\e - \frac{\phi(v_\e)}{\phi'(v_\e)}\right)(\pa_t h - \Delta h)\quad \textrm{in}\;\; Q_T\,.
\end{aligned}
\end{equation}
By \eqref{e54},
\begin{equation}\label{e73}
\phi'(v_\e)=\phi(v_\e)^q = \left(\frac{u_\e}h \right)^q\,.
\end{equation}
>From \eqref{e72}, \eqref{e73} and \eqref{e39} we obtain
\[
\begin{aligned}
& \pa_t(h v_\e) -\Delta(h v_\e) \\ &\,\geq \,  - h^q \left(\frac u{u_\e}\right)^q V + \frac{\phi''(v_\e)}{\phi'(v_\e)} |\nabla v_\e|^2 h + \left(v_\e - \frac{\phi(v_\e)-1}{\phi'(v_\e)}\right)(\pa_t h - \Delta h)\quad \textrm{in}\;\; Q_T\,.
\end{aligned}
\]
In view of \eqref{e38}, \eqref{e37} and \eqref{e12}, the previous inequality implies
\begin{equation}\label{e74}
\pa_t( h v_\e) - \Delta (h v_\e) \geq - h^q \left(\frac{u}{u_\e}\right)^q V\quad \textrm{in}\;\; Q_T\,.
\end{equation}
If $q>0, q\neq 1$, from \eqref{e56} we have that
\[ \phi^{-1}(s) = \frac{s^{1-q} -1}{1-q}\,, \;\; s>0\,,  \]
hence
\begin{equation}\label{e75}
h v_\e = h \phi^{-1}\left(\frac{u_\e}{h} \right)= \frac 1{1-q}(h^q u_\e^{1-q} - h)\quad \textrm{in}\;\; Q_T\,. \
\end{equation}
Let $(x_0, t_0)\in \big[\pa \Omega\times (0, T] \big] \cup \big[\Omega\times \{0\}\big]$. Since $u, h \in C(\bar Q_T)$, in view of \eqref{e39}
we have that
\begin{equation}\label{e76}
u_\e(x_0, t_0)\geq h(x_0, t_0) + \e > h(x_0, t_0)\,.
\end{equation}
>From \eqref{e75} and \eqref{e76} we deduce that
\begin{equation}\label{e77}
\lim_{(x,t)\to (x_0, t_0)} h(x,t) v_\e(x,t) = \frac 1{1-q}\big[h^q(x_0,t_0)u_\e^{1-q}(x_0,t_0) - h(x_0,t_0)\big]\geq 0\,.
\end{equation}
For $q=1$, we have that $\phi^{-1}(s)=\log s$, hence
\begin{equation}\label{e78}
 h v_\e = h \log\left(\frac{u_\e}{h} \right)\quad \textrm{in}\;\; Q_T\,.
\end{equation}
If $h(x_0, t_0)>0$, then we have
\begin{equation}\label{e79}
\lim_{(x,t)\to (x_0, t_0)} h(x,t) v_\e(x,t)=h(x_0, t_0)\log\left(\frac{u_\e(x_0, t_0)}{h(x_0,t_0)}\right)>0\,,
\end{equation}
while if $h(x_0, t_0)=0$, then from \eqref{e78}, since $u_\e\geq \e$, we have that
\begin{equation}\label{e80}
\lim_{(x,t)\to (x_0, t_0)} h(x,t) v_\e(x,t)=0\,.
\end{equation}
>From \eqref{e77}, \eqref{e79} and \eqref{e80} we can infer that $h v_\e\in C^{2,1}(Q_T)\cap C(\bar Q_T),$ and
\begin{equation}\label{e81}
h v_\e \geq 0 \quad \textrm{in}\;\; \big[\pa \Omega\times (0, T] \big] \cup \big[\Omega\times \{0\}\big]\,.
\end{equation}

Note that since
\[ \mathcal S^\Omega\left[ h^q\left(\frac{u}{u_\e}\right)^q |V|\right] \leq \mathcal S^\Omega[ h^q |V|],
\]
we can infer that $\mathcal S^\Omega\left[h^q \left(\frac{u}{u_\e}\right)^q V\right]<\infty$ in $Q_T$ ; furthermore, $h^q\left(\frac{u}{u_\e} \right)^q V\in C(Q_T)$.
Hence, in view of \eqref{e74} and \eqref{e81}, we can apply Proposition \ref{prop1} to obtain
\[ h v_\e \geq - \mathcal S^\Omega\left[h^q \left(\frac{u}{u_\e}\right)^q V\right]\quad \textrm{in}\;\; \, Q_T\,.
\]
Therefore,
\begin{equation}\label{e82}
 v_\e\geq - \frac {1} h \mathcal S^\Omega\left[h^q \left(\frac{u}{u_\e}\right)^q V\right]\quad \textrm{in}\;\; \, Q_T\,.
\end{equation}

\smallskip

We claim that, if $q\geq 1$, then
\begin{equation}\label{e83}
u>0\quad \textrm{in}\;\; Q_T\,.
\end{equation}

\smallskip

In fact, from \eqref{e82} we obtain
\begin{equation}\label{e84}
v_\e \geq - \frac 1 h \mathcal S^\Omega[ h^q V^+]\quad \textrm{in}\;\; Q_T\,.
\end{equation}
Observe that
\[ v_\e = \phi^{-1}\left(\frac{u_\e}{h}\right)\in I_q\,,\;\; - \frac 1 h \mathcal S^\Omega[ h^q V^+]\subset [-\infty, 0]\subseteq \bar I_q\,.  \]
Hence we can apply $\phi$ to both sides of \eqref{e84} to get
\begin{equation}\label{e85}
u_\e \geq h \phi\left(-\frac 1 h \mathcal S^{\Omega}[h^q V^+] \right)\,.
\end{equation}
Letting $\e\to 0^+$ in \eqref{e85} we have
\begin{equation}\label{e86}
u \geq h \phi\left(-\frac 1 h \mathcal S^\Omega[ h^q V^+] \right)\quad \textrm{in}\;\; Q_T\,.
 \end{equation}
Since $\mathcal S^\Omega[h^q V^+](x,t)<\infty$ for every $(x,t)\in Q_T$,  from \eqref{e86} we can infer that \eqref{e83} is satisfied, and the Claim has been shown.

\smallskip

Now, observe that since
\[ v_\e \in I_q, \quad - \frac {1} h \mathcal S^\Omega\left[h^q \left(\frac{u}{u_\e}\right)^q V\right]\in \bar I_q, \]
we can apply $\phi$ to both sides of \eqref{e82} to get
\begin{equation}\label{e87}
u_\e \geq h \phi\left(-\frac 1 h \mathcal S^\Omega\left[h^q \left(\frac{u}{u_\e} \right)^q V \right] \right) \quad \textrm{in}\;\; Q_T\,.
\end{equation}
In view of \eqref{e83}, we have that
\[ \frac{u}{u_\e}\to 1\quad \textrm{in}\;\; Q_T\,\, \textrm{as}\;\; \e\to 0^+\,.\]
Hence, by monotone convergence theorem,
\begin{equation}\label{e88}
\mathcal S^\Omega\left[h^q\left(\frac{u}{u_\e}\right)^q V\right] \to \mathcal S^\Omega\left[h^q V \right] \quad \textrm{in}\;\; Q_T\,\, \textrm{as}\;\, \e\to 0^+\,.
\end{equation}
In particular, we have that
\begin{equation}\label{e89}
-\frac 1 {h(x,t)}\frac{\mathcal S^\Omega[h^q V](x,t)}{h(x,t)}\in \bar I_q\,.
\end{equation}
Letting $\e\to 0^+$ in \eqref{e87} we get
\[u\geq h \phi\left(-\frac 1 h \mathcal S^\Omega[h^q V] \right) \quad \textrm{in}\;\; Q_T\,,\]
from which \eqref{e70} immediately follows. Hence \eqref{e32} and \eqref{e34} have been proved. Furthermore, if $q>1$, from \eqref{e70} we have
\[\phi\left(-\frac 1h \mathcal S^\Omega[h^q V] \right)\leq \frac  u h <\infty\,,
\]
thus
\[ -\frac 1 h \mathcal S^\Omega[h^q V] < \frac 1{q-1}\,,\]
which gives \eqref{e33}\,.

\medskip

Assume that $0<q<1$\,. By the same arguments as in the case $q\geq 1$ we can arrive to \eqref{e82}\,. We can apply $\phi$ to both sides of \eqref{e82} to get
\begin{equation}\label{e90}
u_\e \geq h \phi\left(- \frac 1 h \mathcal S^\Omega\left[h^q\left(\frac{u}{u_\e} \right)^q V \right] \right)\,.
\end{equation}
We have
\[ \frac{u}{u_\e} \to \chi_u \quad \textrm{in}\;\; Q_T\,\,\, \textrm{as}\;\; \e \to 0^+\,.\]
This combined with \eqref{e90} gives
\begin{equation}\label{e92}
u\geq h \phi \left( -\frac 1h \mathcal S^\Omega[\chi_u h^q V ]\right) \quad \textrm{in}\;\; Q_T\,,
\end{equation}
which is equivalent to \eqref{e35}\,.

\medskip

Assume now that $q<0$. For every $\e>0$ we define
\[v_\e := \phi^{-1}\left(\frac{u}{h_\e} \right)\quad \textrm{in}\;\; Q_T\,,\]
where $h_\e:= h+\e\,.$ Since $\frac u{h_\e}>0$ in $Q_T$, we obtain $v_\e\in C^{2,1}(Q_T)$. We extend the function
\begin{equation}\label{e90b}
\phi^{-1}(s)=\frac{s^{1-q}-1}{1-q}, \,\,\, s>0\,,
\end{equation}
by putting $\phi^{-1}(0)=-\frac 1{1-q}\,.$ Since $\frac u{h_\e}\in C(\bar Q_T), \frac{u}{h_\e}\geq 0$ in $\bar Q_T$, we have that $v_\e\in C(\bar Q_T)\,.$

>From \eqref{e40} we have that
\[ u\leq h < h_\e\quad \textrm{in}\;\; \big[ \pa \Omega\times (0, T]\big]\cup \big[\Omega\times\{0\} \big]\,.\]
Hence
\[v_\e \leq \phi^{-1}(1)=0 \quad \textrm{in}\;\; \big[ \pa \Omega\times (0, T]\big]\cup \big[\Omega\times\{0\} \big]\,, \]
therefore,
\begin{equation}\label{e94}
h_\e v_\e \leq 0 \quad \textrm{in}\;\; \big[ \pa \Omega\times (0, T]\big]\cup \big[\Omega\times\{0\} \big]\,.
\end{equation}
In view of \eqref{e40} we have that $u=h_\e \phi(v_\e)$ satisfies
\begin{equation}\label{d1}
\pa_t u - \Delta u + V u^q \leq \pa_t h_\e - \Delta h_\e \quad \textrm{in}\;\; Q_T\,.
\end{equation}
Hence from Lemma \ref{l3} and \eqref{e54} we have that
\begin{equation}\label{d2}
\pa_t (h_\e v_\e) - \Delta (h_\e v_\e)  \leq - h_\e^q V\quad \textrm{in}\;\; Q_T\,.
\end{equation}
Since $q<0$ we have
\[\mathcal S^\Omega[ h_\e^q |V|] \leq \mathcal S^\Omega[ h ^q |V|]\quad \textrm{in}\;\; Q_T\,,\]
so $\mathcal S^\Omega[h_\e ^q V]<\infty$ in $Q_T$\,. Thus, in view of \eqref{d2} and \eqref{e94} we can apply Proposition \ref{prop1a} with $g=-h_\e^q V$ to get
\[h_\e v_\e \leq - \mathcal S^\Omega[h_\e^q V]\quad \textrm{in}\;\; Q_T\,,\]
therefore
\begin{equation}\label{d3}
v_\e \leq -\frac 1{h_\e} \mathcal S^\Omega[h_\e^q V]\quad \textrm{in}\;\; Q_T\,.
\end{equation}
Since $v_\e > -\frac 1 {1-q},$ it follows that
\begin{equation}\label{d4}
-\frac 1{1-q}<-\frac 1{h_\e}\mathcal S^\Omega[h_\e^q V]\leq \infty\,.
\end{equation}
So, we can apply $\phi$ to both sides of \eqref{d3}, and we obtain
\[\phi(v_\e) \leq \phi\left(-\frac 1{h_\e}\mathcal S^\Omega[h_\e^q V] \right)\quad \textrm{in}\;\; Q_T\,,
\]
that is
\[u \leq h_\e \left[1-(1-q)\frac 1{h_\e}\mathcal S^\Omega[h_\e^q V] \right]^{\frac 1{1-q}}\quad \textrm{in}\;\; \, Q_T\,.\]
Therefore,
\begin{equation}\label{d5}
u \leq h_\e \left[1-(1-q)\frac 1{h_\e}\mathcal S^\Omega[h_\e^q V^+]+(1-q)\frac 1{h_\e}\mathcal S^\Omega[h_\e^q V^-] \right]^{\frac 1{1-q}}\,.
\end{equation}
Since $0<h<h_\e$ in $Q_T$ and $q<0$, we have that
\[\frac 1{h_\e} \mathcal S^\Omega[h_\e^q V^-] \leq \frac 1{h} \mathcal S[h^q V^-]\quad \textrm{in}\;\; Q_T\,.  \]
Letting $\e\to 0^+$, by the monotone convergence theorem we obtain
\begin{equation}\label{d6}
\mathcal S^\Omega[ h_\e^q V^+] \to \mathcal S^\Omega[h^q V^+]\quad \textrm{in}\;\; Q_T\,.
\end{equation}
Since $\mathcal S^\Omega[h^q V]$ is well-defined in $Q_T$, letting $\e \to 0^+$ in \eqref{d5},  we have \eqref{e36}.
Since we have assumed that $u>0$ in $Q_T$, from \eqref{e36} it follows \eqref{e33}\,. \hfill $\square$

\section{Proof of Theorems \ref{teo1}, \ref{teo3} and \ref{teo4}}\label{dim}
\setcounter{equation}{0}
\noindent{\it Proof of Theorem \ref{teo1}\,.}  At first, let us show that it is not restrictive to suppose that $f$ is locally  Lipschitz continuous in $Q_T$. In fact, suppose only that $f$ is  continuous in $Q_T$. Let $q>0$. Choose a sequence of nonnegative locally Lipschitz functions $\{f_n \}$ such that
\begin{equation}\label{d8}
f_n\leq f \quad \textrm{in}\;\; \, Q_T\,,
\end{equation}
and
\begin{equation}\label{d7}
f_n\to f \quad \textrm{in}\;\; \; Q_T\;\;\, \textrm{as}\;\;\, n\to \infty\,.
\end{equation}
Set
\begin{equation}\label{d7b}
h_n:= \mathcal R^\Omega[f_n]\,.
\end{equation}
Note that for every $n\in \mathbb N$, $h_n\in C^{2,1}(Q_T)\cap C(\bar Q_T)$ solves \eqref{e38} and \eqref{e37}.
Moreover, we have that
\begin{equation}\label{d7bb}
h_n \leq h,\;\; h_n\to h\quad \textrm{in}\;\; \; Q_T\;\;\, \textrm{as}\;\;\, n\to \infty\,,
\end{equation}
where $h$ is defined in \eqref{e31b}. Since
\[\mathcal S^\Omega[h_n^q |V|]  \leq \mathcal S^\Omega[h^q |V|]\quad \textrm{in}\;\; Q_T\,,\]
we obtain that $\mathcal S^\Omega[h_n^q V]<\infty$ in $Q_T$ for every $n\in \mathbb N$. We have that
\begin{equation}\label{d9}
\mathcal S^\Omega[h_n^q V]\to \mathcal S^\Omega[h^q V]\quad \textrm{in}\;\; Q_T\,,
\end{equation}
and that
\[\mathcal S^\Omega[\chi_u h_n^q V]\to \mathcal S^\Omega[ \chi_uh^q V]\quad \textrm{in}\;\; Q_T\,.\]
In view of \eqref{d8} we deduce that
\begin{equation}\label{d10}
\pa_t u -\Delta u + V u^q \geq f_n\quad \textrm{in}\;\; Q_T\,.
\end{equation}
Therefore, if  \eqref{e32}-\eqref{e35} hold with $h$ replaced by $h_n$ given by \eqref{d7b} and $f$ replaced by $f_n$, then,  thanks to \eqref{d7bb} and \eqref{d9}, we have that \eqref{e32}, \eqref{e34} and \eqref{e35} hold with $h$ given by \eqref{e31b}.
Moreover, we get
\begin{equation}\label{d11}
-(q-1)\mathcal S^\Omega[h^q V]\leq h\quad \textrm{in}\;\; Q_T\,.
\end{equation}
However, from \eqref{e34} it follows that \eqref{d11} must hold with a strict inequality; thus, \eqref{e33} has been shown.

If $q<0$, then the claim follows arguing in the same way, if instead of condition \eqref{d8} we require that
\begin{equation}\label{d8c}
f_n\geq f \quad \textrm{in}\;\; \, Q_T\,.
\end{equation}
Hence, for all $q\neq 0$, we can assume that $f$ is locally Lipschitz continuous in $Q_T$.
Now, let $q>0$. Choose a sequence of subsets $\{\Omega_n\}\subset\subset \Omega$ such that
\begin{equation}\label{d12}
\Omega_n\;\; \textrm{is relatively compact, connected, open and with}\,\,\pa \Omega_n\,\, \textrm{smooth for every}\,\, n\in \mathbb N,
\end{equation}
\begin{equation}\label{d13}
\Omega_n\subset \Omega_{n+1}\;\;\, \textrm{for every}\;\; n\in \mathbb N\,,\,\,\,\cup_{n=1}^\infty \Omega_n\,=\, \Omega\,.
\end{equation}
We have that $h_n:=\mathcal R^{\Omega_n}[f; u_0]\in C^{2,1}(\Omega_n\times (0, T])\cap C(\bar\Omega_n\times[0, T])$, and
\begin{equation}\label{d14}
\left\{
\begin{array}{ll}
\,   \pa_t h_n -\Delta h_n \,=\, f
&\textrm{in}\,\,\Omega_n\times (0, T]
\\& \\
\textrm{ } h_n \, = 0& \textrm{in\ \ } \pa\Omega_n\times (0, T]
\\& \\
\textrm{ } h_n\, = u_0 & \textrm{in\ \ } \Omega_n\times \{0\}\,.
\end{array}
\right.
\end{equation}
We can always take $n$ big enough so that $f\not\equiv 0$ in $\Omega_n$, and so,
\[0 < h_n < \infty\quad \textrm{in}\;\; Q_T\,.\]
By the monotone convergence theorem,
\[  h_n \to h = \mathcal R^\Omega[f; u_0]\quad \textrm{in}\;\; Q_T\,,\,\, \textrm{as}\;\; n\to \infty\,.\]
In view of \eqref{e30} and \eqref{d14} we have that
\begin{equation}\label{d15}
\left\{
\begin{array}{ll}
\,   \pa_t u -\Delta u \,+ V u^q\geq \, \pa_t h_n - \Delta h_n
&\textrm{in}\,\,\Omega_n\times (0, T]
\\& \\
\textrm{ }u \, \geq  h_n & \textrm{in\ \ } \pa\Omega_n\times (0, T]
\\& \\
\textrm{ }u\geq h_n& \textrm{in\ \ } \Omega_n\times \{0\}
\\& \\
\textrm{ }u \, \geq 0& \textrm{in\ \ } \Omega_n\times (0, T]\,.
\end{array}
\right.
\end{equation}
By Theorem \ref{teo2},
\begin{equation}\label{d16}
u\geq \left\{
\begin{array}{ll}
\,  h_n e^{-\frac 1{h_n}\mathcal S^{\Omega_n}[h_n V]}
&\textrm{if}\,\,q=1\,,
\\& \\
\textrm{ }h_n \left\{1+(q-1)\frac 1{h_n}\mathcal S^{\Omega_n}[h_n^q V] \right\}^{-\frac 1{q-1}} & \textrm{if\ \ } q>1
\\& \\
\textrm{ } h_n \left\{1+(q-1)\frac 1{h_n}\mathcal S^{\Omega_n}[\chi_n h_n^q V] \right\}^{-\frac 1{q-1}}_+ & \textrm{if\ \ } 0<q<1\,
\end{array}
\right.
\end{equation}
in $\Omega_n\times (0, T]$, where $\chi_n:= \chi_u|_{\Omega_n}\,.$ Moreover,
\begin{equation}\label{d17}
1+(q-1)\frac 1{h_n}\mathcal S^{\Omega_n}[h_n^q V]>0\,.
\end{equation}
By the monotone convergence theorem,
\[ \mathcal S^{\Omega_n}[h_n^q V^{\pm}]\to \mathcal S^{\Omega}[h^q V^{\pm}]\quad \textrm{in}\;\; Q_T\;\; \textrm{as}\;\; n\to \infty\,,\]
and
\[ \mathcal S^{\Omega_n}[\chi_n h_n^q V^{\pm}]\to \mathcal S^{\Omega}[\chi_u h^q V^{\pm}]\quad \textrm{in}\;\; Q_T\;\; \textrm{as}\;\; n\to \infty\,.\]
Passing to the limit as $n\to \infty$ in \eqref{d16} gives \eqref{e32}, \eqref{e34} and \eqref{e35}\,.
Let $q>1$. Then from \eqref{d17} we have that
\[1+(q-1)\frac 1{h}\mathcal S^{\Omega}[h^q V]\geq 0\,.
\]
However, since $-\frac 1{q-1}<0$ and $\frac u h<\infty$, the previous inequality yields \eqref{e33}\,.

\smallskip
It remains to prove \eqref{e35}\,. Let $q<0$. Note that since $f$ is locally Lipschitz in $Q_T$, $\mathcal R^{\Omega}[f]\in C^{2,1}(Q_T)$. In fact, for every relatively compact subset $\Omega'\subset\Omega$ with $\pa \Omega'$ smooth, we clearly have that $\mathcal R^{\Omega'}[f]\in C^{2,1}(\Omega'\times (0, T])$. Moreover, the function $w:=\mathcal R^{\Omega}[f]-\mathcal R^{\Omega'}[f]$ solves in the weak sense
\begin{equation}\label{d45}
\pa_t w - \Delta w = 0 \quad \textrm{in}\;\; \Omega'\times (0, T]\,.
\end{equation}
Hence, by standard regularity results, $w\in C^{2,1}(\Omega'\times (0, T])$. Therefore, $\mathcal R^{\Omega}[f]\in C^{2,1}(\Omega'\times (0, T])\,.$ Since $\Omega'$ was arbitrary, the claim follows. For any $\e>0$ define
$$
h_\e := \e + \mathcal R^{\Omega}[f; u_0]\,.
$$
We have that
\[\pa_t h_\e - \Delta h_\e \,=\, f\quad \textrm{in}\;\; Q_T\,.\]
Since $u>0, h_\e>0$ in $Q_T$, the function $v_\e:= \phi^{-1}\left(\frac{u}{h_\e} \right)\in C^{2,1}(Q_T)$. By the same arguments as in the proof of Theorem \ref{teo2}, we obtain
\begin{equation}\label{h33}
\pa_t( h_\e v_\e) - \Delta (h_\e v_\e) \leq - h_\e^q V  \quad \textrm{in}\;\; Q_T\,.
\end{equation}
>From \eqref{e90b} we get
\begin{equation}\label{d18}
h_\e v_\e = h_\e \phi^{-1}\left(\frac u{h_\e} \right) = h_\e^q \frac{u^{1-q}- h_\e^{1-q}}{1-q}\,.
\end{equation}
Observe that
\begin{equation}\label{m11}
u=0 \quad \textrm{in}\;\; \pa \Omega\times (0, T]\,,
\end{equation}
and
\begin{equation}\label{m12}
u(x,0)\leq u_0(x) \quad \textrm{for all}\;\; x\in \Omega\,.
\end{equation}
Moreover,
 \begin{equation}\label{m13}
h_\varepsilon>\e \quad \textrm{in}\;\; \pa \Omega\times (0, T]\,,
\end{equation}
and
\begin{equation}\label{m14}
h_\e(x,0)=\e + u_0 \quad \textrm{for all}\;\; x\in \Omega\,.
\end{equation}
>From \eqref{d18}, \eqref{m11}-\eqref{m14} we can infer that
\begin{equation}\label{d19}
h_\e v_\e \leq 0 \quad \textrm{in}\;\; \big[\pa \Omega\times (0, T]\big]\cup\big[\Omega\times\{0\}\big]\,.
\end{equation}
Moreover, from \eqref{e49} and fact that $h_\varepsilon>\varepsilon$ it follows that
\begin{equation}\label{d20}
\lim_{x\to \pa_\infty M} \sup_{t\in (0, T]} h_\e(x,t) v_\e(x,t)=0\,.
\end{equation}
Therefore, we can apply Proposition \ref{prop1a} with $g=-h_\e^q V$ to get
\begin{equation}\label{h36}
h_\e v_\e \leq  - \mathcal S^{\Omega}[h_\e^q V]\quad \textrm{in}\;\; Q_T\,.
\end{equation}
Letting $\e\to 0^+$, the thesis follows by the same arguments as in the proof of Theorem \ref{teo2}-$(iv)$. This completes the proof. \hfill $\square$

\medskip

\noindent{\it Proof of Theorem \ref{teo3}\,.}  Let $\{\Omega_n\}$ be a sequence of domains as in \eqref{d12}-\eqref{d13}. Let $q\geq 1$. For every $n\in \mathbb N$, let $h_n\in C^{2,1}(Q_T)\cap C(\bar Q_T)$ be the solution of problem
\[
\left\{
\begin{array}{ll}
\,   \pa_t h_n -\Delta h_n \,=\,0
&\textrm{in}\,\,\Omega_n\times (0, T]
\\& \\
\textrm{ } h_n \, = u& \textrm{in\ \ } \pa\Omega_n\times (0, T]
\\& \\
\textrm{ } h_n\, = u & \textrm{in\ \ } \Omega_n\times \{0\}\,.
\end{array}
\right.
\]
In view of \eqref{e43} and \eqref{e45}, by the maximum principle,
\[ h_n >0 \quad \textrm{in}\;\; Q_T\,.\]
Thanks to \eqref{e83}, we can infer that $u(x)>0$ for all $x\in \Omega_n, t\in (0, T]$; therefore,
$u(x)>0$ for all $(x,t)\in Q_T$\,.

\smallskip

Let $q=1$. Set $h\equiv 1, v:= \log u.$ As in the proof of Theorem \ref{teo2}, we have
\[
\pa_t v - \Delta v \geq - V\quad \textrm{in}\;\; Q_T\,.
\]
>From \eqref{e43} we can deduce that
\[v\geq 0 \quad \textrm{in}\;\; \big[\pa\Omega \times (0, T]\big]\cup \big[\Omega\times\{0\}\big]\,,\]
and
\[\liminf_{x\to \pa_\infty M}\inf_{t\in (0, T]} v(x,t) \geq 0\,.\]
Thus, we can apply Proposition \ref{prop1} with $g=-V$, and we have
\begin{equation}\label{d22}
\log u(x,t) = v(x,t) \geq -\mathcal S^{\Omega}[V](x,t)\quad \textrm{for all}\;\; (x,t)\in Q_T\,.
\end{equation}
>From \eqref{d22}, inequality \eqref{e44} immediately follows.

\smallskip

Now, let $q>1$. Set
$$\a_n:= \inf_{[\pa \Omega_n\times (0, T]]\cup [\Omega_n\times \{0\}]} u\,.$$
In view of \eqref{e45} we have that
\begin{equation}\label{d23}
\lim_{n\to \infty} \a_n\,=\, \infty\,.
\end{equation}
We can apply Theorem \ref{teo2} with $h\equiv \a_n$. Therefore,
\begin{equation}\label{d24}
\begin{aligned}
 u \geq &  \a_n \left\{ 1 + (q-1) \a_n^{q-1}\mathcal S^{\Omega_n}[V] \right\}^{-\frac 1{q-1}} \\
=& \{ \a_n^{-(q-1)}+ (q-1) \mathcal S^{\Omega_n}[V]\}^{-\frac 1{q-1}}\quad \textrm{in}\;\; \Omega_n\times (0, T]\,,
\end{aligned}
\end{equation}
and
\begin{equation}\label{d25}
-(q-1) \mathcal S^{\Omega_n}[V] < \a_n^{-(q-1)}\quad \textrm{in}\;\; \Omega_n\times (0, T]\,.
\end{equation}
Hence, letting $n\to \infty$ in \eqref{d25} we get $\mathcal S^\Omega[V](x)\geq 0.$ Therefore, by the monotone convergence theorem, \eqref{d24} implies \eqref{e47}. Since $u(x)<\infty$, \eqref{e46} follows.

\smallskip

Now, let $0<q<1$. We set
\[\phi(v):=[(1-q) v]_+^{\frac 1{1-q}}\,,\quad v\in \mathbb R.\]
Thus
\[ \phi'(v)>0,\;\; \phi''(v) >0 \quad \textrm{for all}\;\; v>0\,.\]
Moreover, \eqref{e54} holds.
Consider a sequence $\{\e_n\}\subset (0,\infty)$ with $\e_n\to 0$ as $n\to \infty$. For every $n\in \mathbb N$ define
\[ u_n:= u + \e_n,\;\; v_n:=\phi^{-1}(u_n)\,.\]
In view of Remark \ref{ossl3} with $h\equiv 1$, by the same arguments as in the proof of Theorem \ref{teo2}, we have
\[\pa_t v_n - \Delta v_n \geq -\left(\frac{u_n}{u}\right)^q V\quad \textrm{in}\;\; \Omega_n\times (0, T]\,.
\]
Since
$$
v_n > 0\quad \textrm{in}\;\; \big[\pa\Omega_n\times (0, T]\big]\cup \big[\Omega_n\times\{0\}\big]\,,
$$
by Proposition \ref{prop1},
\begin{equation}\label{d30}
v_n \geq -\mathcal S^{\Omega_n}\left[\left(\frac{u_n}{u}\right)^q V \right]\quad \textrm{in}\;\; \Omega_n\times (0, T]\,.
\end{equation}
Letting $n\to \infty$, by the monotone convergence theorem we get
\[\phi^{-1}(u) \geq - \mathcal S[\chi_{u} V]\quad \textrm{in}\;\; Q_T\,,\]
which is equivalent to \eqref{e48}.

\smallskip

Now, let $q<0$. For every $n\in \mathbb N$ set
\[\nu_n:= \sup_{\big[\Omega_n\times\{0\} \big] \cup \big[\pa\Omega_n\times (0, T] \big]} u\,\,. \]
In view of \eqref{e49a} and \eqref{e49} we have that
\begin{equation}\label{h20}
\lim_{n\to \infty} \nu_n\,=\,0\,.
\end{equation}
We can apply Theorem \ref{teo2} in $\Omega_n$ with $h\equiv \nu_n$ to obtain
\begin{equation}\label{h21}
u(x,t)\leq \big\{\nu_n^{1-q} - (1-q) \mathcal S^{\Omega_n}[V](x,t) \big\}^{\frac 1{1-q}} \quad \textrm{for all}\;\; (x,t)\in Q_T\,.
\end{equation}
Letting $n\to \infty$ in \eqref{h21} we get \eqref{e51}. Moreover, since $u>0$ in $Q_T$, we obtain \eqref{e50h}.
This completes the proof. \hfill $\square$

\medskip
\bigskip

In order to prove Theorem \ref{teo4} we use the standard method of sub-- and supersolutions; namely, if there exists $\underline u, \overline u\in C^{2,1}(Q_T)\cap C(\bar Q_T)$ such that
\begin{equation}\label{d36a}
0\leq \underline u \leq \overline u\quad \textrm{in}\;\; \, Q_T,
\end{equation}
\begin{equation}\label{d37b}
\underline u = 0, \quad \overline u \geq 0 \quad \textrm{in}\;\; \pa \Omega\times (0, T]\,,
\end{equation}
\begin{equation}\label{m5}
\underline u \leq u_0 \leq \overline u \quad \textrm{in}\;\; \Omega\times\{0\}\,.
\end{equation}
and
\begin{equation}\label{d38c}
\pa_t \underline u - \Delta \underline u + V \underline u^q \leq f\, \quad \textrm{in}\;\; Q_T\,,
\end{equation}
\begin{equation}\label{d38b}
\pa_t \overline u -\Delta \overline u + V \overline u^q \geq f \, \quad \textrm{in}\;\; Q_T\,,
\end{equation}
then there exists a solution $u\in C^{2,1}(Q_T)\cap C(\bar Q_T)$ of  problem \eqref{d39}
such that
\begin{equation}\label{d39b}
\underline u \leq u \leq \overline u \quad \textrm{in}\;\; Q_T\,.
\end{equation}

\noindent{\it Proof of Theorem \ref{teo4}\,.} We limit ourselves to prove the statement $(ii)$, since the statement $(i)$ can be proved in a similar and simpler way.

Let
\[\overline u\equiv h = \mathcal R^\Omega[f; u_0]\,.\]
In view of the regularity assumptions on $f$ and on $\partial\Omega$, we have that $\overline u\in C^{2,1}(Q_T)\cap C(\bar Q_T)$ solves
\[\left\{
\begin{array}{ll}
\,   \pa_t \overline u -\Delta \overline u \, =\, f
&\textrm{in}\,\,Q_T
\\& \\
\textrm{ }\overline u \, = 0& \textrm{in\ \ } \pa\Omega\times (0, T] \,.
\\& \\
\textrm{ }\overline u \, = u_0& \textrm{in\ \ } \Omega\times \{0\}\,.
\end{array}
\right.
\]
Moreover, since $V\geq 0, f\geq 0$, we have that $\overline u$ satisfies \eqref{d38b}. Hence $\overline u$ is a supersolution of problem \eqref{d39}.

Now, we look for a subsolution $\underline u$ of problem \eqref{d39}. To this aim, define
\[ \underline u := h - \lambda^q \mathcal S^\Omega[ h^q V]\quad \textrm{in}\;\; Q_T\,,\]
where $\lambda>0$ is a positive parameter to be fixed in the sequel. Thanks to \eqref{d37} we have that if we take
\begin{equation}\label{d40}
0<\lambda<-\frac{q(1-q)^{\frac 1 q}}{1-q}\,,
\end{equation}
then
$$\underline u>0\quad \textrm{in}\;\; Q_T\,.$$ Hence, \eqref{d36a} holds.
We claim that $\underline u \in C^{2,1}(Q_T)\cap C(\bar Q_T)\,.$ In fact, for every relatively compact subset $\Omega'\subset\Omega$ with $\pa \Omega'$ smooth, since $h>0$ in $\bar \Omega'$, we have that $\mathcal S^{\Omega'}[h^q V]\in C^{2,1}(\Omega'\times (0, T])$. Moreover, the function $w:=\mathcal S^{\Omega}[h^q V]-\mathcal S^{\Omega'}[h^q V]$ solves
\eqref{d45} in the weak sense. Hence, by standard regularity results, $w\in C^{2,1}(\Omega'\times (0, T])$. Therefore, $\mathcal S^{\Omega}[h^q V]\in C^{2,1}(\Omega'\times (0, T])\,.$ Since $\Omega'$ was arbitrary, the claim follows. Furthermore, since $h\in C(\bar Q_T)$ and $h=0$ in $\big[\pa\Omega\times (0, T] \big]\cup \big[\Omega\times\{0\}\big]$, using \eqref{d37} we can deduce that $\mathcal S^\Omega[h^q V]\in C(\bar Q_T)$ and $\mathcal S^{\Omega}[h^q V]=0\,$ in $\big[\pa\Omega\times (0, T] \big]\cup \big[\Omega\times\{0\}\big]$.

Now, let us show that $\underline u$ satisfies \eqref{d38c}. Note that
$$
\pa_t \underline u -\Delta \underline u +  V \underline u^q = f - \lambda^q h^q V + \underline u^q V\quad \textrm{in}\,\, Q_T.
$$
Hence, since $V\geq 0$ and $q<0$, \eqref{d38} follows, if we show that
\[\lambda h \leq \underline u, \]
that is
\begin{equation}\label{d47}
\mathcal S^\Omega[ h^q V] \leq \lambda^{-q}(1-\lambda) h\,.
\end{equation}
Now, it is easily checked that \eqref{d37} yields \eqref{d47}, by taking $\lambda=\frac 1{1-\frac 1 q}$. Consequently, there exists a solution $u\in C^{2,1}(Q_T)\cap C(\bar Q_T)$ of problem \eqref{d39} such that
\eqref{d39b} is satisfied. Therefore,
\[ u\geq \underline u = h -\lambda^q \mathcal S^\Omega[h^q V] = h - (1-\frac 1 q)^{-q}\mathcal S^\Omega[h^q V]\geq \frac 1{1-\frac 1 q} h \quad \textrm{in}\;\; Q_T\,.\]
This combined with Theorem \ref{teo1}-$(iv)$ gives \eqref{d38}. The proof is complete. \hfill $\square$

\section{Proof of Theorems \ref{teo5} and \ref{teo6}}\label{dim2}
\noindent{\it Proof of Theorem \ref{teo5}\,.} By the same arguments as in the proof of Theorem \ref{teo1}, and using the same notations, we can infer that, for any $\e>0$,
\eqref{h33} and \eqref{d19} hold. In view of \eqref{d18} and \eqref{h30} we have that for any $\e>0$
\begin{equation}\label{h34}
\limsup_{x\to \pa_\infty M} \frac{\sup_{t\in (0, T]} h_\e(x,t) v_\e(x,t)}{|Z(x)|} \leq 0\,.
\end{equation}
Due to \eqref{h34} we can apply Proposition \ref{prop1r} with $g=-h_\e^q V$ to deduce \eqref{h36}. Thus the conclusion follows as in the proof of Theorem \ref{teo1}. \hfill $\square$

\medskip

\noindent{\it Proof of Theorem \ref{teo6}\,.}
Choose a sequence of not relatively compact domains $\{\Omega_n\}_{n\in \mathbb N}$ with smooth boundary such that
\[ \Omega_{n}\subset \Omega_{n+1},\; \bar \Omega_n\subset \Omega\, \quad \textrm{for every}\,\, n\in \mathbb N,\;\; \cup_{n=1}^\infty \Omega_n\,=\, \Omega\,.\]

For every $n\in \mathbb N$ set
\begin{equation}\label{k1}
\nu_n:= \sup_{\big[\Omega_n\times\{0\} \big] \cup \big[\pa\Omega_n\times (0, T] \big]} u\,\,.
\end{equation}
In view of \eqref{e49a} we have that
\begin{equation}\label{h20}
\lim_{n\to \infty} \nu_n\,=\,0\,.
\end{equation}
For each $n\in \mathbb N$ set $h:= \nu_n$. Since $u>0, h>0$ in $Q_T$, the function $v:= \phi^{-1}\left(\frac{u}{h} \right)\in C^{2,1}(Q_T)$; here $\phi^{-1}$ is given by \eqref{e90b}. By the same arguments as in the proof of Theorem \ref{teo2}, we obtain
\begin{equation}\label{k33}
\pa_t( h v) - \Delta (h v) \leq - h^q V  \quad \textrm{in}\;\; Q_T\,.
\end{equation}
>From \eqref{e90b} we get
\begin{equation}\label{k34}
h v = h \phi^{-1}\left(\frac u{h} \right) = h^q \frac{u^{1-q}- h^{1-q}}{1-q}\,.
\end{equation}
>From \eqref{k1} we can infer that
\begin{equation}\label{k35}
h v \leq 0 \quad \textrm{in}\;\; \big[\pa \Omega_n\times (0, T]\big]\cup\big[\Omega_n\times\{0\}\big]\,.
\end{equation}
Moreover, due to \eqref{k34} and \eqref{h32} we have that
\begin{equation}\label{k36}
\limsup_{x\to \pa_\infty M} \frac{\sup_{t\in (0, T]} h(x,t)v(x,t) }{|Z(x)|}\, \leq\, 0\,.
\end{equation}
Therefore, for each $n\in \mathbb N$ we can apply can apply Proposition \ref{prop1r}  with $g=-h^q V$ to get
\begin{equation}\label{k37}
h v \leq  - \mathcal S^{\Omega}[h^q V]\quad \textrm{in}\;\; \Omega_n\times (0, T]\,.
\end{equation}
Hence by Theorem \ref{teo2} in $\Omega_n$ with $h\equiv \nu_n$ we obtain
\begin{equation}\label{k38}
u(x,t)\leq \big\{\nu_n^{1-q} - (1-q) \mathcal S^{\Omega_n}[V](x,t) \big\}^{\frac 1{1-q}} \quad \textrm{for all}\;\; (x,t)\in \Omega_n\times (0, T]\,.
\end{equation}
Letting $n\to \infty$ in \eqref{k38}, using \eqref{k1}, we get \eqref{e51}. Moreover, since $u>0$ in $Q_T$, we obtain \eqref{e50h}.
This completes the proof.


\begin{thebibliography}{999}

\bibitem{Aubin} T. Aubin, "Some nonlinear problems in Riemannian Geometry", Springer (1998)\,.

%\bibitem{AB} D.G. Aronson, P. Besala, {\it Uniqueness of solutions to the Cauchy problem for parabolic equations}, J. Math.
%Anal. Appl. {\bf 13} (1966), 516--526\,.

\bibitem{BPT} C. Bandle, M. A. Pozio, A. Tesei, {\it The Fujita exponent for the Cauchy problem in the hyperbolic space}, J. Diff. Eq. {\bf 251} (2011), 2143--2163 .

\bibitem{BC} H. Brezis, X. Cabré, {\it Some simple nonlinear PDE's without solutions}, Boll. Unione Mat. Ital. {\bf 8}, Ser. 1-B (1998) 223-262\,.

\bibitem{BK} H. Brezis, S. Kamin, {\it Sublinear elliptic equations in $\mathbb R^n$,} Manuscr. Math. {\bf 74} (1992), 87--106\,.

\bibitem{FrV} M. Frazier, I.E. Verbitsky, {\it Global Green's function estimates}, Around the Reseach of Vladimir Maz'ya III, Analysis and Applications, Ed. Ari Laptev Math. Series {\bf 13}, Springer, (2010), 105--152\,.


\bibitem{Grig} \newblock A. Grigoryan, \emph{ Analytic and geometric background of
recurrence and non-explosion of the Brownian motion on Riemannian
manifolds}, \newblock Bull. Amer. Math. Soc. {\bf 36} (1999),
135--249.

\bibitem{Grig3} A. Grigor'yan, ``Heat Kernel and Analysis on Manifolds'', AMS/IP Studies in Advanced Mathematics, 47, American Mathematical Society, Providence, RI; International Press, Boston, MA, 2009.

\bibitem{GrigH2} A. Grigor'yan, W. Hansen, {\it Lower estimates for a perturbed Green function}, J. Anal. Math. {\bf 104} (2008) , 25--58\,.

\bibitem{GrigH} A. Grigor'yan, W. Hansen, {\it Lower estimates for perturbed Dirichlet solutions}, preprint (1999)\,.

\bibitem{GrigV} A. Grigor'yan, I. E. Verbitsky, {\it Pointwise estimates of solutions to semilinear elliptic equations and inequalities}, J. d'Anal. Math. (to appear)\,.

%\bibitem{IKO} A. M. Il'in, A. S. Kalashnikov, O. A. Oleinik, {\it Linear equations of the second order of parabolic
%type}, Russian Math. Surveys {\bf 17} (1962), 1--144\,.


\bibitem{KV} N.J. Kalton, I.E. Verbitsky, {\it Nonlinear equations and weighted norm inequalities}, Trans. Amer. Math. Soc. {\bf 351} (1999), 3441--3497\,.

\bibitem{MMP} P. Mastrolia, D. Monticelli, F. Punzo, {\it Nonexistence of solutions to parabolic differential inequalities with a potential on Riemannian manifolds}, Math. Ann., DOI 10.1007/s00208-016-1393-2\,.

%\bibitem{Pinch} Y. Pinchover, {\it On uniqueness and nonuniqueness of the positive
%Cauchy problem for parabolic equations with unbounded
%coefficients}, Math.Z. {\bf 223} (1996), 569--586\,.
%
%
%\bibitem{Ti} A. N. Tihonov, {\it
%Th\'eor\`emes d'unicit\'e pour l'\'equation de la chaleur}, Mat.
%Sb. {\bf 42} (1935), 199--215\,.

\bibitem{Pu1} F. Punzo, {\it Blow-up of solutions to semilinear parabolic equations on Riemannian manifolds with negative sectional
curvature}, J. Math. Anal. Appl., {\bf 387} (2012), 815--827 .

\bibitem{Zhang} Q. S. Zhang, {\it Blow-up results for nonlinear parabolic equations on manifolds}, Duke Math. J. {\bf 97} (1999), 515--539 .

\end{thebibliography}
\end{document}